\documentclass[11pt]{amsart}

\usepackage{graphicx}
\usepackage{amssymb}
\usepackage{amsfonts}
\usepackage[]{amsmath}
\usepackage[]{epsfig}
\usepackage[]{pstricks}
\usepackage[]{float}
\newpsobject{malla}{psgrid}{subgriddiv=1,griddots=10,gridlabels=6pt}

\newtheorem{theorem}{Theorem}[section]
\newtheorem{lemma}{Lemma}[section]

\newtheorem{remark}{Remark}[section]
\newtheorem{proposition}{Proposition}[section]

\newtheorem{claim}{Claim}[section]
\newtheorem{question}{Question}[section]
\numberwithin{equation}{section}

\newcommand{\R}{\mathbb R}
\newcommand{\Z}{\mathbb Z}
\newcommand{\T}{\mathbb T}

\newcommand{\SL}{\mathcal{L}}
\newcommand{\SM}{\mathcal{M}}
\newcommand{\SN}{\mathcal{N}}
\newcommand{\Ss}{\mathcal{S}}
\newcommand{\ST}{\mathcal{T}}
\newcommand{\SU}{\mathcal{U}}
\newcommand{\ov}{\overline}

\newcommand{\Int}{\int \! \! \! \! \int}
\newcommand{\la}{\langle}
\newcommand{\ra}{\rangle}

\newcommand{\p}{\partial}

\begin{document}

\title[Rough solutions for the periodic NLS-KdV system]{Rough solutions for the 
periodic Schr\"odinger - Kortweg-deVries system}

\author{A. Arbieto}
\address{Alexander Arbieto 
\newline IMPA, Estrada Dona Castorina 110, Rio de Janeiro, 22460--320, 
Brazil.} \email{alexande@impa.br}
\author{A. J. Corcho}
\address{Ad\'an J. Corcho
\newline Universidade Federal de Alagoas, Departamento de Matem\'atica, 
Campus A. C. Sim\~oes, Ta\-buleiro dos Martins, Macei\'o-AL, 
57072-900, Brazil.} \email{adan@mat.ufal.br}
\author{C. Matheus}
\address{Carlos Matheus 
\newline IMPA, Estrada Dona Castorina 110, Rio de Janeiro, 22460--320,
Brazil.} \email{matheus@impa.br}


\date{October 11, 2005.}

\begin{abstract}We prove two new mixed sharp bilinear estimates of 
Schr\"odinger-Airy type. In particular, we obtain the local 
well-posedness of the Cauchy problem of the Schr\"odinger - Kortweg-deVries 
(NLS-KdV) system in the \emph{periodic setting}. Our lowest regularity is 
$H^{1/4}\times L^2$, which is somewhat far from the naturally expected 
endpoint $L^2\times H^{-1/2}$. This is a novel phenomena related to the 
periodicity condition. Indeed, in the continuous case, Corcho and Linares 
proved local well-posedness for the natural endpoint 
$L^2\times H^{-\frac{3}{4}+}$. 

Nevertheless, we conclude the global well-posedness of the NLS-KdV 
system in the energy space $H^1\times H^1$ using our local 
well-posedness result and three conservation laws discovered by M. 
Tsutsumi.
\end{abstract}

\maketitle


\section{Introduction}\label{s.intro}

The interaction of a short-wave $u=u(x,t)$ and a long-wave $v=v(x,t)$ in 
fluid mechanics (and plasma physics) is governed by the Schr\"odinger - 
Kortweg-de Vries (NLS-KdV) system

\begin{equation}\label{e.nls-kdv}
\begin{cases}
i\partial_tu + \partial_x^2u = \alpha uv + \beta |u|^2u,& t\in {\mathbb R},\\
\partial_tv + \partial_x^3v + \tfrac{1}{2}\partial_x(v^2) = 
\gamma \partial_x(|u|^2),\\
u(x,0)=u_0(x),\;v(x,0)=v_0(x),
\end{cases} \end{equation}
where $u=u(x,t)$ is a complex-valued function, $v=v(x,t)$ is a 
real-valued function and $\alpha,\;\beta,\;\gamma$\; are real 
constants.\footnote{The case $\beta = 0$ of the NLS-KdV system 
occurs in the study of the resonant interaction between short and 
long capillary-gravity waves on water channels of uniform finite 
depth and in a diatomic lattice system. For more details about 
these physical applications, see~\cite{Funakoshi}, 
\cite{Nishikawa}, \cite{Kawahara} and \cite{Yajima-Satsuma}.}

This motivates the study of the local and global well-posedness of 
the Cauchy problem for the NLS-KdV system with rough initial 
data.\footnote{Benilov and Burtsev in~\cite{Benilov} showed that 
the NLS-KdV is not completely integrable. In particular, the 
solvability of~(\ref{e.nls-kdv}) depends on the theory of 
non-linear dispersive equations.} 

The central theme of this paper is the local and global 
well-posedness theory of the NLS-KdV system in the periodic 
setting (i.e., $x\in\T$); but, in order to motivate our subsequent 
results, we recall some known theorems in the non-periodic 
setting. 

In the continuous context (i.e., $x\in\R$), Corcho and 
Linares~\cite{Corcho} showed the local well-posedness of the 
NLS-KdV for initial data $(u_0,v_0)\in H^k(\R)\times H^s(\R)$ with 
$k\geq 0$, $s>-3/4$ provided that $k-1\leq s\leq 2k-\frac{1}{2}$ 
for $k\leq 1/2$ and $k-1\leq s < k+\frac{1}{2}$ for $k>1/2$. It is 
worth to point out that the lowest regularity obtained by Corcho 
and Linares is $k=0$ and $s=-\frac{3}{4} +$. In the non-resonant 
case $\beta\neq 0$, it is reasonable to expect that the NLS-KdV is 
locally well-posed in $L^2\times H^{-\frac{3}{4}+}$: the nonlinear 
Schr\"odinger (NLS) equation with cubic term $(|u|^2 u)$ is 
globally well-posed in $H^s(\R)$ for $s\geq 0$ and ill-posed below 
$L^2(\R)$; similarly, the Kortweg-de Vries (KdV) equation is 
globally well-posed in $H^s(\R)$ for $s>-3/4$ and ill-posed in 
$H^s(\R)$ for $-1\leq s < -3/4$. Also, using three conserved 
quatinties for the NLS-KdV flow, M. Tsutsumi~\cite{MTsutsumi} 
showed global well-posedness for initial data $(u_0,v_0)\in 
H^{s+\frac{1}{2}}(\R)\times H^s(\R)$ with $s\in\Z_+$ and Corcho, 
Linares~\cite{Corcho}, assuming $\alpha\gamma>0$, showed global 
well-posedness in the energy space $H^1(\R)\times 
H^1(\R)$.\footnote{Pecher~\cite{Pecher} announced the global 
well-posedness of the NLS-KdV system (with $\alpha\gamma>0$) in 
the continuous setting for initial data $(u_0,v_0)\in 
H^s(\R)\times H^s(\R)$, for $3/5<s<1$ in the resonant case 
$\beta=0$ and $2/3<s<1$ in the non-resonant case $\beta\neq 0$. 
The proof is based on two refined bilinear estimates and the 
I-method of Colliander, Keel, Staffilani, Takaoka and Tao.} 
 
The point of view adopted by Corcho and Linares in order to prove 
their local well-posedness result is to use a basic strategy to 
treat, in both continuous and periodic contexts, the 
low-regularity study of dispersive equations (such as NLS and 
KdV): one considers the Fourier restriction norm method introduced 
by Bourgain in~\cite{Bourgain}; then, they showed two new mixed 
bilinear estimates for the coupling terms of the NLS-KdV system 
(namely, $uv$ and $\partial_x(|u|^2)$) in certain Bourgain's 
spaces, which implies that an equivalent integral equation can be 
solved by Picard's fixed point method (in other words, the 
operator associated to the integral equation is a contraction in 
certain Bourgain spaces). Coming back to the periodic setting, 
before stating our results, we advance that, although our efforts 
are to obtain similar well-posedness theorems, the periodic case 
is more subtle than the continuous context: since the cubic NLS is 
globally well-posed (resp., ill-posed) in $H^s(\T)$ for $s\geq 0$ 
(resp. $s<0$) and the KdV is globally well-posed (resp., 
ill-posed) in $H^s(\T)$ for $s\geq -1/2$ (resp., $s<-1/2$), it is 
reasonable to expect $L^2(\T)\times H^{-1/2}(\T)$ as the lowest 
regularity for the local well-posedness results; but, surprisingly 
enough, the endpoint of the bilinear estimates for the coupling 
terms $uv$, $\p_x(|u|^2)$ in the periodic setting is $(k,s)=(1/4, 
0)$, i.e., our lowest regularity is $H^{1/4}\times L^2$ (see the 
propositions~\ref{p.uv},~\ref{p.du2}, theorem~\ref{t.A} and 
remark~\ref{r.1} below). We refer the reader to the 
section~\ref{s.remarks} for a more detailed comparasion between 
the well-posedness results for the NLS-KdV system in the periodic 
and non-periodic settings (as well as a couple of questions 
motivated by this discussion). 
 
Now, we introduce some notations. Let $U(t) = e^{it\p_x^2}$ and 
$V(t) = e^{-t\p_x^3}$ be the unitary groups associated to the 
linear Schr\"odinger and the Airy equations, respectively. Given 
$k,s,b\in\R$, we define the spaces $X^{k,b}$ and $Y^{s,b}$ via the 
norms 
\begin{equation*}
\begin{split}
\|f\|_{X^{k,b}} &:= \left(\sum\limits_{n\in\Z} \la n \ra^{2k} \la 
\tau+n^2 
\ra^{2b} |\widehat{f}(n,\tau)|^2\right)^{1/2} \\ 
&= \|U(-t) f\|_{H_t^b(\R,H_x^k)} 
\end{split}
\end{equation*}
\begin{equation*}
\begin{split}
\|g\|_{Y^{s,b}} &:= \left(\sum\limits_{n\in\Z} \la n \ra^{2k} 
\la\tau-n^3\ra^{2b} |\widehat{g}(n,\tau)|^2\right)^{1/2} \\
&= \|V(-t) g\|_{H_t^n(\R,H_x^s)}
\end{split}
\end{equation*}
where $\la\cdot\ra:= 1+|\cdot|$ and $\widehat{f}$ is the Fourier 
transform of $f$ in both variables $x$ and $t$:
\begin{equation*}
\widehat{f}(n,\tau) = (2\pi)^{-1}\int_{\R\times\T} e^{-it\tau}
e^{-ixn} f(x,t) dtdx
\end{equation*}
and, given a time interval $I$, we define $X^{k,b}(I)$ and 
$Y^{s,b}(I)$ via the (restriction in time) norms
\begin{equation*}
\|f\|_{X^{k,b}(I)} = \inf\limits_{\widetilde{f}|_I = f}
\|\widetilde{f}\|_{X^{k,b}} \quad \textrm{and} \quad 
\|g\|_{Y^{s,b}(I)} = \inf\limits_{\widetilde{g}|_I = g}
\|\widetilde{g}\|_{Y^{s,b}}
\end{equation*}

The study of the periodic dispersive equations (e.g, KdV) has been 
based around iteration in the Bourgain spaces (e.g., $Y^{s,b}$) 
with $b=1/2$. Since we are interested in the continuity of the 
flow associated to the NLS-KdV system and the Bourgain spaces with 
$b=1/2$ do not control the $L_t^{\infty}H_x^s$, we consider the 
slightly smaller spaces $\widetilde{X}^{k}$, $\widetilde{Y}^s$ 
defined by the norms
\begin{equation*}
\|u\|_{\widetilde{X}^{k}}:= \|u\|_{X^{k,1/2}} + \|\la n\ra^k
\widehat{u}(n,\tau)\|_{L_n^2 L_{\tau}^1} \quad \textrm{and} \quad 
\|v\|_{\widetilde{Y}^{s}}:= \|v\|_{Y^{s,1/2}} + \|\la n\ra^s
\widehat{v}(n,\tau)\|_{L_n^2 L_{\tau}^1} 
\end{equation*}
and, given a time interval $I$, we define the spaces 
$\widetilde{X}^{k}(I)$, $\widetilde{Y}^s(I)$ via the restriction 
in time norms
\begin{equation*}
\|f\|_{\widetilde{X}^k(I)} = \inf\limits_{\widetilde{f}|_I = f}
\|\widetilde{f}\|_{\widetilde{X}^k} \quad \textrm{and} \quad 
\|g\|_{Y^{s,b}(I)} = \inf\limits_{\widetilde{g}|_I = g}
\|\widetilde{g}\|_{\widetilde{Y}^s}
\end{equation*}

Also, we introduce the companion spaces $Z^k$ and $W^s$ via the 
norms
\begin{equation*}
\|u\|_{Z^k}:= \|u\|_{X^{k,-1/2}} + \left\|\frac{\la n\ra^k
\widehat{u}(n,\tau)}{\la\tau+n^2\ra}\right\|_{L_n^2 L_{\tau}^1} 
\quad \textrm{and} \quad \|v\|_{W^{s}}:= \|v\|_{Y^{s,-1/2}} + 
\left\|\frac{\la n\ra^s
\widehat{v}(n,\tau)}{\la\tau-n^3\ra}\right\|_{L_n^2 L_{\tau}^1}
\end{equation*}

Denote by $\psi$ a non-negative smooth bump function supported in
$[-2,2]$ with $\psi = 1$ on $[-1,1]$ and $\psi_{\delta}(t):=\psi(t/ \delta)$ 
for any $\delta>0$. Also, let $a\pm$ be a number slightly larger (resp., 
smaller) than $a$.  At this point, we are ready to state our main results. The 
fundamental technical propositions are the following two sharp bilinear for the 
coupling terms of the NLS-KdV system:  

\begin{proposition}\label{p.uv}For any $s\geq 0$ and $k-s\leq 3/2$, 
\begin{equation}\label{e.uv} \|uv\|_{Z^k}\lesssim 
\|u\|_{X^{k,\frac{1}{2}-}}\|v\|_{Y^{s,\frac{1}{2}}} + 
\|u\|_{X^{k,\frac{1}{2}}}\|v\|_{Y^{s,\frac{1}{2}-}}. 
\end{equation}  
Furthermore, the estimate~(\ref{e.uv}) fails if either $s<0$ or $k-s>3/2$. 
More precisely, if the bilinear estimate $\|uv\|_{X^{k,b-1}}\lesssim 
\|u\|_{X^{k,b}} \|v\|_{Y^{s,b}}$ with $b=1/2$ holds then $s\geq 0$ and 
$k-s\leq 3/2$. 
\end{proposition}  
\begin{proposition}\label{p.du2}For any $k>0$, $1+s\leq 4k$ and $-1/2\leq k-s$, 
\begin{equation}\label{e.du2} \|\p_x(u_1 \ov{u_2})\|_{W^s}\lesssim 
\|u_1\|_{X^{k,\frac{1}{2}-}} \|u_2\|_{X^{k,\frac{1}{2}}} + 
\|u\|_{X^{k,\frac{1}{2}}}\|v\|_{X^{k,\frac{1}{2}-}}. 
\end{equation}  
Furthermore, the estimate~(\ref{e.du2}) fails if either $1+s>4k$ or $k-s<-1/2$. 
More precisely, if the bilinear estimate $\|\p_x(u_1 \ov{u_2})\|_{Y^{s,-1/2}} 
\lesssim \|u_1\|_{X^{k,1/2}}\|u_2\|_{X^{k,1/2}}$ holds then $1+s\leq 4k$ and 
$-1/2\leq k-s$. 
\end{proposition}   

Using these bilinear estimates for the coupling terms $uv$ and $\p_x(|u|^2)$, 
we show the main theorem of this paper, namely, we prove the following local 
well-posedness result:  

\begin{theorem}\label{t.A}The periodic NLS-KdV~(\ref{e.nls-kdv}) is locally 
well-posed in $H^k(\T)\times H^s(\T)$ whenever $s\geq 0$, $-1/2\leq k-s\leq 
3/2$ and $1+s\leq 4k$. I.e., for any $(u_0,v_0)\in H^k(\T)\times H^s(\T)$, 
there exists a positive time $T=T(\|u_0\|_{H^k},\|v_0\|_{H^s})$ and a unique 
solution $(u(t),v(t))$ of the NLS-KdV system~(\ref{e.nls-kdv}) satisfying  
\begin{equation*} 
(\psi_T(t) u, \psi_T(t) v)\in \widetilde{X}^k\times\widetilde{Y}^s, 
\end{equation*}  
\begin{equation*} 
(u,v)\in C([0,T],H^k(\T)\times H^s(\T)). 
\end{equation*}  
Moreover, the map $(u_0,v_0)\mapsto (u(t),v(t))$ is locally Lipschitz from 
$H^k(\T)\times H^s(\T)$ into $C([0,T],H^k(\T)\times H^s(\T))$, whenever 
$k,s\geq 0$, $-1/2\leq k-s\leq 3/2$ and $1+s\leq 4k$. 
\end{theorem}  

\begin{remark}\label{r.1}As we pointed out before, the endpoint of our sharp 
bilinear estimates and, consequently, our local well-posedness result is 
$H^{1/4}\times L^2$. Since the endpoint of the sharp well-posedness theory for 
the periodic NLS is $L^2$ and for the periodic KdV is $H^{-1/2}$, we are 
somewhat far from the naturally expected endpoint $L^2\times H^{-1/2}$ for the 
local in time theory for the NLS-KdV system (although, our bilinear estimates 
are optimal). This leads us to ask about possible ill-posedness results in this 
gap between $H^{1/4}\times L^2$ and $L^2\times H^{-1/2}$. For precise 
statements and some comparision with the continuous setting, see the 
section~\ref{s.remarks}. 
\end{remark}   

\begin{remark}\label{r.2} It is easy to see that the NLS-KdV 
system~(\ref{e.nls-kdv}) system is ill-posed for $k<0$. Indeed, if we put 
\begin{equation*} 
\begin{cases} 
u:=e^{-it}w,\\ v\equiv \alpha^{-1} \in H^s(\mathbb{T}), \forall 
s\in \R, 
\end{cases} 
\end{equation*} 
the system (\ref{e.nls-kdv}) becomes into the equation 
\begin{equation*} 
\begin{cases} 
iw_t+\partial_x^2w=\beta |w|^2w,\\ 
\partial_x (|w|^2)=0,\\ 
w_0(x)=u_0\in H^k(\mathbb{T}), 
\end{cases} 
\end{equation*} 
which is not locally-well posed (ill-posed) below $L^2({\mathbb T})$ in the 
sense that the data-solution map is not uniformly continuous. 
\end{remark}  

Using this local well-posedness result and three conserved quantities for the 
NLS-KdV flow, it will be not difficult to prove the following global 
well-posedness theorem in the energy space $H^1(\T)\times H^1(\T)$:  

\begin{theorem}\label{t.B}Let $\alpha,\beta,\gamma\in\R$ be such that 
$\alpha\gamma>0$ and $(u_0,v_0)\in H^1(\T)\times H^1(\T)$. Then, the unique 
solution in the theorem~\ref{t.A} can be extended to the time interval $[0,T]$ 
for any $T>0$. 
\end{theorem}  

To close this introduction, we give the outline of the paper. In 
section~\ref{s.examples} we give counter-examples for the bilinear estimates of 
the coupling terms, when the indices $k$ and $s$ satisfies (at least) one of 
the following inequalities: $s<0$, $k-s>3/2$, $1+s>4k$ or $k-s<-1/2$. In 
section~\ref{s.bilinear} we complete the proof of the propositions~\ref{p.uv} 
and~\ref{p.du2} by establishing the claimed bilinear estimates for the terms 
$uv$ and $\p_x(|u|^2)$. In section~\ref{s.thmA} we use propositions~\ref{p.uv} 
and~\ref{p.du2} to show that the integral operator associated to the NLS-KdV 
system is a contraction in the space $\widetilde{X}^{k}([0,T])\times 
\widetilde{Y}^{s}([0,T])$ (for sufficiently small $T>0$) when $k,s\geq 0$, 
$-1/2\leq k-s\leq 3/2$ and $1+s\leq 4k$. In particular, we obtain the desired 
local well-posedness statement in theorem~\ref{t.A}. In section~\ref{s.thmB} 
we make a standard use of three conserved quantities for the NLS-KdV flow to 
obtain the global well-posedness result of theorem~\ref{t.B} in the energy 
space $H^1(\T)\times H^1(\T)$. In section~\ref{s.remarks} we make some 
questions related to the gap between the expected $L^2\times H^{-1/2}$ 
endpoint regularity and our lowest regularity $H^{1/4}\times L^2$ for the 
local well-posedness for the periodic NLS-KdV system; also, we compare the 
known theorems in the continuous setting with the periodic setting. Finally, 
in the appendix, we collect some standard facts about linear and multilinear 
estimates associated to the cubic NLS and the KdV equations (which were used 
in the proof of theorem~\ref{t.A}) and we show that the NLS-KdV flow preserves 
three quantities controlling the $H^1$ norms of $u(t)$ and $v(t)$ (this is the 
heart of the proof of theorem~\ref{t.B}).   

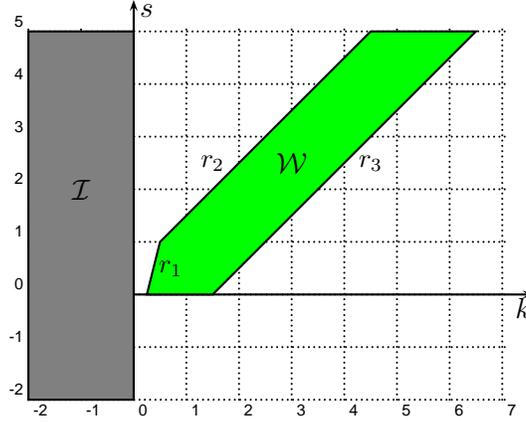
\begin{figure}[ht] 
\centering 
\psset{unit=0.7cm} 
\begin{pspicture}(-2,-2)(7,5)%
\malla%
\psline{->}(-2,0)(7.6,0)%
\psline{->}(0,-2)(0,5.6)%
\pspolygon[fillstyle=solid,fillcolor=green,linewidth=0.8pt] 
(0.25,0)(1.5,0)(6.5,5)(4.5,5)(0.5,1)(0.25,0)%
\rput(3,2.5){$\mathcal{W}$}%
\pspolygon[fillstyle=solid,fillcolor=gray,linewidth=0.8pt] 
(0,-2)(-2,-2)(-2,5)(0,5)(0,-2)%
\rput(-1,2){$\mathcal{I}$}%
\rput(7.4,-0.25){$k$}%
\rput(0.25,5.4){$s$}%
\rput(4.5,2.5){\small{$r_3$}}%
\rput(1.5,2.5){\small{$r_2$}}%
\rput(0.7,0.5){\small{$r_1$}}%
\end{pspicture} 
\vspace{0.5cm} 
\caption{Well-posedness results for periodic NLS-KdV system. The region 
$\mathcal{W}$, limited for the lines $r_1: s=4k-1$, $r_2: s=k+\frac{1}{2}$ and 
$r_3: s=k-\frac{3}{2}$, contains indices $(k,s)$ for which local 
well-posedness is achieved in Theorem \ref{t.A}. The region 
$\mathcal{I}$ show the ill-posedness results commented in Remark \ref{r.2}.} 
\end{figure}    

\section{Counter-Examples}\label{s.examples} We start with some 
counter-examples for the bilinear estimate in proposition~\ref{p.uv} when 
$s<0$ or $k-s>3/2$:  

\begin{lemma}\label{l.example-uv}$\|uv\|_{X^{k,b-1}}\leq 
\|u\|_{X^{k,b}}\cdot\|v\|_{Y^{s,b}}$ (with $b=1/2$) implies $s\geq 0$ and 
$k-s\leq 3/2$. 
\end{lemma}  

\begin{proof}Fix $N\gg 1$ a large integer. Firstly, 
we show that $\|uv\|_{X^{k,b-1}}\leq \|u\|_{X^{k,b}}\cdot\|v\|_{Y^{s,b}}$ 
(with $b=1/2$) implies $s\geq 0$. Define 
\begin{displaymath} b_n= \left\{ \begin{array}{ll} 1 & \textrm{if $n=N$}\\ 
0 & \textrm{otherwise} \end{array} \right. 
\end{displaymath} 
and 
\begin{displaymath} a_n= \left\{ \begin{array}{ll} 1 & \textrm{if 
$n=\frac{-N^2 -N}{2}$}\\ 
0 & \textrm{otherwise} \end{array} \right. 
\end{displaymath}  

Let $u$ and $v$ be defined by $\widehat{u}(n,\tau)=a_n\chi_1(\tau+n^2)$ and 
$\widehat{v}(n,\tau)=b_n\chi_1(\tau-n^3)$, where $\chi_1$ is the characteristic 
function of the interval $[-1,1]$.  Now let's go to the calculations. By 
definition of the Bourgain space $X^{k,b}$, 
$$\|u v\|_{X^{k,b-1}}=\left\|\frac{\la n \ra^k}{\la \tau+n^2 \ra^{1/2}} 
\widehat{u}\ast\widehat{v}\right\|_{L^2_{n,\tau}}.$$  Hence, 
$$\|u v\|_{X^{k,b-1}}= \left\|\frac{\la n \ra^k}{\la \tau+n^2 \ra^{1/2}} 
\sum\limits_{n_1}\int d\tau_1 \ a_{n-n_1} \ \chi_1((\tau-\tau_1)+(n-n_1)^2) \ 
b_{n_1}\ \chi_1(\tau_1-n_1^3)\right\|_{L^2_{n,\tau}}.$$  

Recall the following numerical expression:  
\begin{equation}\label{e.1} 
\left(\tau_1-n_1^3\right)+\left((\tau-\tau_1)+(n-n_1)^2\right) 
-\left(\tau+n^2\right) = -n_1^3 + n_1^2 - 2 n n_1. 
\end{equation}  

Taking into account that $b_{n_1}\neq 0$ iff $n_1=N$, $a_{n-n_1}\neq 0$ 
iff $n=\frac{-N^2 +N}{2}$, $\chi_1(\tau_1-n_1^3)\neq 0$ 
iff $|\tau_1-n_1^3|\leq 1$ and $\chi_1((\tau-\tau_1)+(n-n_1)^2)\neq 0$ 
iff $|(\tau-\tau_1)+(n-n_1)^2|\leq 1$, we conclude, from a direct substitution 
of these data into~(\ref{e.1}), that  
\begin{equation}\label{e.example-uv} 
\|u v\|_{X^{k,b-1}}\approx N^{2k}. 
\end{equation}  

On the other hand, it is not difficult to see that 
 
\begin{equation}\label{e.example-u} 
\|u\|_{X^{k,b}}= \|\la n \ra^k \la \tau+n^2 \ra^{1/2} \ a_n \ 
\chi_1(\tau+n^2)\|_{L^2_{n,\tau}} \approx N^{2k}, 
\end{equation} 
and 
\begin{equation}\label{e.example-v} 
\|v\|_{Y^{s,b}} = \|\la n \ra^s \la \tau-n^3 \ra^{1/2} b_n \ 
\chi_1(\tau-n^3)\|_{L^2_{n,\tau}} \approx N^s. 
\end{equation} 

Putting together the 
equations~(\ref{e.example-uv}),~(\ref{e.example-u}),~(\ref{e.example-v}), 
we obtain that the bilinear estimate implies 

$$N^{2k}\lesssim N^{2k}\cdot N^s,$$ 
which is possible only if $s\geq 0$. 

Secondly, we prove that $\|uv\|_{X^{k,b-1}}\leq 
\|u\|_{X^{k,b}}\cdot\|v\|_{Y^{s,b}}$ (with $b=1/2$) implies 
$k-s\leq 3/2$. 
 
Define 
\begin{displaymath} 
b_n= \left\{ \begin{array}{ll} 
1 & \textrm{if $n=N$}\\ 
0 & \textrm{otherwise} 
\end{array} \right. 
\end{displaymath} 
and 
\begin{displaymath} 
a_n= \left\{ \begin{array}{ll} 
1 & \textrm{if $n=0$}\\ 
0 & \textrm{otherwise} 
\end{array} \right. 
\end{displaymath} 

Let $u$ and $v$ be defined by 
$\widehat{u}(n,\tau)=a_n\chi_1(\tau+n^2)$ and 
$\widehat{v}(n,\tau)=b_n\chi_1(\tau-n^3)$, where $\chi_1$ is the 
characteristic function of the interval $[-1,1]$.

Using the definitions of the Bourgain $X^{k,b}$ and $Y^{s,b}$ 
spaces and the algebraic relation~(\ref{e.1}), we have 
$$\|u v\|_{X^{k,b-1}}\approx \frac{N^k}{N^{3/2}},$$ 
$$\|u\|_{X^{k,b}}\approx 1,$$ 
and 
$$\|v\|_{Y^{s,b}}\approx N^s.$$ 

Hence, the bilinear estimate says 
$$N^k\lesssim N^s N^{3/2},$$ 
which is only possible if $k-s\leq 3/2$. 
\end{proof}

We consider now some counter-examples for the bilinear estimate in 
proposition~\ref{p.du2} when $1+s> 4k$ or $k-s<-1/2$: 

\begin{lemma}\label{l.example-du2}$\|\p_x (u_1 \ov{u_2})\|_{Y^{s,b-1}}\leq 
\|u_1\|_{X^{k,b}}\cdot\|u_2\|_{X^{k,b}}$ (with $b=1/2$) implies 
$1+s\leq 4k$ and $k-s\geq -1/2$. 
\end{lemma}

\begin{proof}Fix $N\gg 1$ a large integer. 
Firstly, we prove that $\|\p_x (u_1 \ov{u_2})\|_{Y^{s,b-1}}\leq 
\|u_1\|_{X^{k,b}}\cdot\|u_2\|_{X^{k,b}}$ (with $b=1/2$) implies 
$1+s\leq 4k$. Define 
\begin{displaymath} 
b_n= \left\{ \begin{array}{ll} 
1 & \textrm{if $n=\frac{-N^2-N}{2}$}\\ 
0 & \textrm{otherwise} 
\end{array} \right. 
\end{displaymath} 
and 
\begin{displaymath} 
a_n= \left\{ \begin{array}{ll} 
1 & \textrm{if $n=\frac{-N^2 +N}{2}$}\\ 
0 & \textrm{otherwise} 
\end{array} \right. 
\end{displaymath} 

Let $u_1$ and $u_2$ be defined by 
$\widehat{u_1}(n,\tau)=a_n\chi_1(\tau+n^2)$ and 
$\widehat{u_2}(n,\tau)=b_n\chi_1(\tau+n^2)$, where $\chi_1$ is the 
characteristic function of the interval $[-1,1]$. 
 
By definition of the Bourgain space $Y^{s,b}$, 
$$\|\p_x(u_1\ov{u_2})\|_{Y^{s,b-1}}= 
\left\|\frac{\la n \ra^s} {\la \tau-n^3 \ra^{1/2}} \; n \; 
(\widehat{u_1}\ast\widehat{\ov{u_2}})\right\|_{L^2_{n,\tau}}.$$ 
 
Hence, if one uses that 
$\widehat{\ov{u}}(n,\tau)=\ov{\widehat{u}(-n,-\tau)}$, it is not 
difficult to see that 
$$\|\p_x(u_1\ov{u_2})\|_{Y^{s,b-1}}= 
\left\|\frac{ n \; \la n \ra^s}{\la \tau-n^3 \ra^{1/2}} 
\sum\limits_{n_1}\int d\tau_1 \ a_{n-n_1} \ 
\chi_1((\tau-\tau_1)+(n-n_1)^2) \ b_{-n_1}\ 
\chi_1(-\tau_1+n_1^2)\right\|_{L^2_{n,\tau}}.$$ 
  
Note the following numerical expression: 
   
\begin{equation}\label{e.2} 
\left(\tau-n^3\right)-\left((\tau-\tau_1)+(n-n_1)^2\right) 
+\left(-\tau_1+n_1^2\right) = -n^3 - n^2 + 2 n_1 n. 
\end{equation} 
    
Taking into account that $b_{-n_1}\neq 0$ iff 
$n_1=\frac{N^2+N}{2}$, $a_{n-n_1}\neq 0$ iff $n=N$, 
$\chi_1(-\tau_1+n_1^2)\neq 0$ iff $|-\tau_1+n_1^2|\leq 1$ and 
$\chi_1((\tau-\tau_1)+(n-n_1)^2)\neq 0$ iff 
$|(\tau-\tau_1)+(n-n_1)^2|\leq 1$, we conclude, from a direct 
substitution of these data into~(\ref{e.2}), that 
     
\begin{equation}\label{e.example-du2} 
\|\p_x(u_1\ov{u_2})\|_{Y^{s,b-1}}\approx N^{1+s}. 
\end{equation} 
      
On the other hand, it is not difficult to see that 
       
\begin{equation}\label{e.u_1} 
\|u_1\|_{X^{k,b}}= \|\la n \ra^k \la \tau+n^2 \ra^{1/2} \ a_n \ 
\chi_1(\tau+n^2)\|_{L^2_{n,\tau}} \approx N^{2k}, 
\end{equation} 
and 
\begin{equation}\label{e.u_2} 
\|u_2\|_{X^{k,b}} = \|\la n \ra^k \la \tau+n^2 \ra^{1/2} b_n \ 
\chi_1(\tau+n^2)\|_{L^2_{n,\tau}} \approx N^{2k}. 
\end{equation} 
        
Putting together the 
equations~(\ref{e.example-du2}),~(\ref{e.u_1}),~(\ref{e.u_2}), we 
obtain that the bilinear estimate implies 
$$N^{1+s}\lesssim N^{2k}\cdot N^{2k},$$ 
which is possible only if $1+s\leq 4k$. 
	 
Secondly, we obtain that $\|\p_x (u_1 \ov{u_2})\|_{Y^{s,b-1}}\leq 
\|u_1\|_{X^{k,b}}\cdot\|u_2\|_{X^{k,b}}$ (with $b=1/2$) implies 
$k-s\geq -1/2$. 
	  
Define 
\begin{displaymath} 
b_n= \left\{ \begin{array}{ll} 
1 & \textrm{if $n=-N$}\\ 
0 & \textrm{otherwise} 
\end{array} \right. 
\end{displaymath} 
and 
\begin{displaymath} 
a_n= \left\{ \begin{array}{ll} 
1 & \textrm{if $n=0$}\\ 
0 & \textrm{otherwise} 
\end{array} \right. 
\end{displaymath}  

Let $u_1$ and $u_2$ be defined by 
$\widehat{u_1}(n,\tau)=a_n\chi_1(\tau+n^2)$ and 
$\widehat{u_2}(n,\tau)=b_n\chi_1(\tau+n^2)$, where $\chi_1$ is the 
characteristic function of the interval $[-1,1]$. 
 
Using the definitions of the Bourgain $X^{k,b}$ and $Y^{s,b}$ 
spaces and the algebraic relation~(\ref{e.1}), we have 
$$\|\p_x (u_1\ov{u_2})\|_{Y^{s,b-1}}\approx \frac{N^{1+s}}{N^{3/2}},$$ 
$$\|u_1\|_{X^{k,b}}\approx 1,$$ 
and 
$$\|u_2\|_{Y^{s,b}}\approx N^k.$$  

Hence, the bilinear estimate says 
$$N^{1+s}\lesssim N^k N^{3/2},$$ 
which is only possible if $k-s\geq -1/2$. 
\end{proof}  

\section{Bilinear Estimates for the Coupling Terms}\label{s.bilinear} 
This section is devoted to the proof of our basic tools, that is, 
the sharp bilinear estimates~\ref{p.uv},~\ref{p.du2} for the 
coupling terms of the NLS-KdV system. We begin by showing some 
elementary calculus lemmas; next, using Plancherel and duality, 
the claimed bilinear estimates reduce to controlling some weighted 
convolution integrals, which is quite easy from these lemmatas. 
 
\subsection{Preliminaries} 
The first elementary calculus lemma to be used later is: 
 
\begin{lemma}\label{l.calculus-1} 
\begin{equation*} 
\int_{-\infty}^{+\infty}\frac{d\kappa}{\la\kappa\ra^{\theta} 
\la\kappa-a\ra^{\widetilde{\theta}}} \lesssim 
\frac{\log (1+\la a \ra)}{\la a \ra^{\theta+\widetilde{\theta}-1}} 
\end{equation*} 
where $\theta,\widetilde{\theta}>0$ and 
$\theta+\widetilde{\theta}>1$ 
\end{lemma}  

\begin{proof}Clearly we can assume that $|a|\gg 1$. In this case, we divide the
domain of integration into the regions $I_1:=\{|\kappa|\ll |a|\}$, 
$I_2:=\{|\kappa|\sim |a|\}$ and $I_3:=\{|\kappa|\gg |a|\}$. Since $\kappa\in 
I_1$ implies $\la \kappa-a\ra\gtrsim \la a\ra\geq \la \kappa\ra$, 
$\kappa\in I_2$ implies $\la \kappa\ra\sim\la a\ra$ and
$x\in I_3$ implies $\la \kappa-a\ra\gtrsim \la \kappa\ra$, we obtain 
\begin{equation*}
\begin{split}
\int_{-\infty}^{+\infty}\frac{d\kappa}{\la\kappa\ra^{\theta} 
\la\kappa-a\ra^{\widetilde{\theta}}} &= \int_{I_1}
\frac{d\kappa}{\la\kappa\ra^{\theta} 
\la\kappa-a\ra^{\widetilde{\theta}}} + 
\int_{I_2}
\frac{d\kappa}{\la\kappa\ra^{\theta} 
\la\kappa-a\ra^{\widetilde{\theta}}} + 
\int_{I_3}
\frac{d\kappa}{\la\kappa\ra^{\theta} 
\la\kappa-a\ra^{\widetilde{\theta}}}\\ 
&\lesssim \frac{1}{\la a\ra^{\theta+\widetilde{\theta}-1}}
\int_{I_1}\frac{d\kappa}{\la\kappa\ra} + 
\frac{1}{\la a\ra^{\theta+\widetilde{\theta}-1}}
\int_{I_2}\frac{d\kappa}{\la\kappa-a\ra} + 
\int_{I_3}\frac{d\kappa}{\la\kappa-a\ra^{\theta+\widetilde{\theta}}} \\ 
&\lesssim \frac{\log (1+\la a \ra)}{\la a \ra^{\theta+\widetilde{\theta}-1}}.
\end{split}
\end{equation*} 
\end{proof} 
 
The second lemma is a well-known fact concerning the convergence 
of series whose terms are the values of certain polynomials along 
the integer numbers:\footnote{This lemma is essentially contained 
in the work~\cite{KPV2} of Kenig, Ponce and Vega on bilinear 
estimates related to the KdV equation.}  

\begin{lemma}\label{l.1.1}For any constant $\theta > 1/3$,  
$$ 
\sum\limits_{m\in\Z}\frac{1}{\la p(m) \ra^{\theta}} \leq 
C(\theta)<\infty, 
$$ 
where $p(x)$ is a cubic polynomial of the form $p(x):= 
x^3+ex^2+fx+g$ with $e,f,g\in\R$. 
\end{lemma}  

\begin{proof} 
We start the proof of the lemma~\ref{l.1.1} with two simple 
observations: defining 
 
\begin{equation*} 
\mathcal{E}:=\{m\in\Z: \ |m-\alpha|\geq 2, \ |m-\beta|\geq 2 
\text{ and } |m-\gamma|\geq 2\} 
\end{equation*} 
and 
\begin{equation*} 
\mathcal{F}:=\Z-\mathcal{E}, 
\end{equation*} 
then 
$$\#\mathcal{F}\leq 12$$ and $$\la (m-\alpha)(m-\beta)(m-\gamma) \ra 
\gtrsim \la m-\alpha \ra \la m-\beta \ra \la m-\gamma \ra$$ for 
any $m\in\mathcal{E}.$ 

In particular, writing $p(x)=(x-\alpha)(x-\beta)(x-\gamma)$, we 
can estimate 
\begin{equation*} 
\begin{split} 
\sum\limits_{m}\frac{1}{p(m)^{\theta}}&\leq 
\sum\limits_{m\in\mathcal{F}}\frac{1}{p(m)^{\theta}} + 
\sum\limits_{m\in\mathcal{E}}\frac{1}{p(m)^{\theta}} \\ 
&\leq 12 + \sum\limits_{m\in\mathcal{E}}\frac{1}{p(m)^{\theta}} \\ 
&\lesssim 12 + \sum\limits_{m\in\mathcal{E}}\frac{1}{\la m-\alpha 
\ra^{\theta} \la m-\beta \ra^{\theta} \la m-\gamma \ra^{\theta}} 
\end{split} 
\end{equation*} 

Now, by H\"older inequality 
\begin{equation*} 
\sum\limits_{m}\frac{1}{p(m)^{\theta}}\lesssim 12 + 
\left(\sum\limits_{m}\frac{1}{\la m-\alpha 
\ra^{3\theta}}\right)^{1/3} \left(\sum\limits_{m}\frac{1}{\la 
m-\beta \ra^{3\theta}}\right)^{1/3} 
\left(\sum\limits_{m}\frac{1}{\la m-\gamma 
\ra^{3\theta}}\right)^{1/3} 
\end{equation*} 

So, the hypothesis $3\theta>1$ implies 
\begin{equation*} 
\sum\limits_{m}\frac{1}{p(m)^{\theta}}\leq C(\theta)<\infty. 
\end{equation*} 
This completes the proof of the lemma~\ref{l.1.1}. 
\end{proof}  

Finally, the third lemma is a modification of the previous one for 
linear polynomials with large coefficients: 

\begin{lemma}\label{l.1.2}For any constant $\theta>1/2$, whenever 
$n_1\in\Z-\{0\}$, $|n_1|\gg 1$, 
\begin{equation*} 
\sum\limits_{n\in\Z; \ |n|\sim |n_1|}\frac{1}{q(n)^{\theta}}\leq 
C(\theta)<\infty, 
\end{equation*} 
where $q(x):=2 n_1 x-n_1^2+r$ with $r\in\R$. 
\end{lemma}  

\begin{proof}The strategy of the proof is the same as before, but since now 
the polynomial $q$ is linear, we have to take a little bit of 
care. The idea is: although the polynomial $q$ has degree $1$, the 
fact that $|n_1|\sim |n|$ means morally that $q$ has degree $2$ in 
this range. So, the exponent of $n$ in the summand is morally 
$2\theta>1$ and, in particular, the series is convergent. This 
intuition can be formalized as follows: we write 
$q(x):=r-n_1^2+ 2 n_1 x = 2 n_1 (x+\delta)$, where 
$\delta = (r-n_1^2)/(2 n_1)$ (of course the assumption $n_1\neq 0$ 
enters here). If we define 
\begin{equation*} 
\mathcal{G}:=\{n\in\Z: \ |n+\delta|\geq 2\} 
\end{equation*} 
and 
\begin{equation*} 
\mathcal{H}:=\Z-\mathcal{G}, 
\end{equation*} 
then 
$$\#\mathcal{H}\leq 4$$ and $$\la 2 n_1(n+\delta) \ra 
\gtrsim \la n \ra \la n-\delta \ra$$ for any $n\in\mathcal{G}$, 
since $|n_1|\sim |n|$. 
 
In particular, we can estimate 
\begin{equation*} 
\begin{split} 
\sum\limits_{|n|\sim |n_1|}\frac{1}{q(n_1)^{\theta}}&\leq 
\sum\limits_{n\in\mathcal{H}}\frac{1}{q(n_1)^{\theta}} + 
\sum\limits_{n\in\mathcal{G},|n|\sim |n_1|}\frac{1}{q(n_1)^{\theta}} \\ 
&\leq 4 + 
\sum\limits_{n\in\mathcal{G},|n|\sim |n_1|}\frac{1}{q(n_1)^{\theta}} \\ 
&\lesssim 4 + \sum\limits_{n\in\mathcal{G},|n|\sim 
|n_1|}\frac{1}{\la n \ra^{\theta} \la n+\delta \ra^{\theta}} 
\end{split} 
\end{equation*} 

Now, by H\"older inequality 
\begin{equation*} 
\sum\limits_{|n|\sim |n_1|}\frac{1}{q(n_1)^{\theta}}\lesssim 4 + 
\left(\sum\limits_{n}\frac{1}{\la n \ra^{2\theta}}\right)^{1/2} 
\left(\sum\limits_{n}\frac{1}{\la n+\delta 
\ra^{2\theta}}\right)^{1/2} 
\end{equation*} 

So, the hypothesis $2\theta>1$ implies 
\begin{equation*} 
\sum\limits_{|n|\sim |n_1|}\frac{1}{q(n)^{\theta}}\leq 
C(\theta)<\infty. 
\end{equation*} 
This completes the proof of the lemma~\ref{l.1.2}. 
\end{proof} 

\subsection{Proof of the proposition~\ref{p.uv}: 
bilinear estimates for the coupling term $uv$} In view of the 
lemma~\ref{l.example-uv}, it suffices to show the bilinear 
estimate: 
 
\begin{lemma}\label{l.sharp-uv} $\|uv\|_{Z^k}\lesssim 
\|u\|_{X^{k,\frac{1}{2}-}}\|v\|_{Y^{s,\frac{1}{2}}} + 
\|u\|_{X^{k,\frac{1}{2}}}\|v\|_{Y^{s,\frac{1}{2}-}}$ whenever 
$s\geq 0$ and $k-s\leq 3/2$. 
\end{lemma}  

\begin{proof} From the definition of $Z^k$, we must show that 
\begin{equation}\label{e.sharp-uv-1} 
\|uv\|_{X^{k,-1/2}}\lesssim 
\|u\|_{X^{k,\frac{1}{2}-}}\|v\|_{Y^{s,\frac{1}{2}}} + 
\|u\|_{X^{k,\frac{1}{2}}}\|v\|_{Y^{s,\frac{1}{2}-}} 
\end{equation} 
and 
\begin{equation}\label{e.sharp-uv-2} 
\left\|\frac{\la n \ra^k \widehat{uv}}{\la \tau+n^2 \ra} 
\right\|_{L^2_n L^1_{\tau}} \lesssim 
\|u\|_{X^{k,\frac{1}{2}-}}\|v\|_{Y^{s,\frac{1}{2}}} + 
\|u\|_{X^{k,\frac{1}{2}}}\|v\|_{Y^{s,\frac{1}{2}-}} 
\end{equation}  

We begin with the estimate~(\ref{e.sharp-uv-1}). By the definition 
of Bourgain's space, 
\begin{equation*}\begin{split} 
\|uv\|_{X^{k,-a}}&= \|\la \tau + {n}^2 \ra^{-a}\la n \ra^k 
\widehat {uv} (n,\tau )\|_{L^2_{\tau}L^2_{n}}\\ 
&=\Bigl \| \frac{\la n \ra^k}{\la \tau + {n}^2 \ra^a} 
\widehat {u}*\widehat{v}(n,\tau )\Bigl \|_{L^2_{\tau}L^2_{n}}\\ 
\end{split} 
\end{equation*} 
 
Let 
\begin{equation*} 
f(\tau,n)=\la \tau + {n}^2\ra^b\la n \ra^k\widehat {u}(n,\tau 
)\quad \text{and}\quad g(\tau,n)=\la \tau - {n}^3\ra^c\la n 
\ra^s\widehat {v}(n,\tau ). 
\end{equation*}  

In particular, by duality, we obtain 
\begin{equation}\label{e.sharp-uv-1-duality} 
\begin{split} 
\|uv\|_{X^{k,-a}}&=\sup\limits_{\|\varphi\|_{L^2_{n,\tau}}\leq 1} 
\sum\limits_{n\in\Z}\int d\tau \frac{\la n \ra^k}{\la \tau + n^2 
\ra^{a}}\bar{\varphi}(n,\tau) \left( \frac{f}{\la\tau+n^2\ra^b \la 
n\ra^k}*\frac{g}{\la\tau-n^3\ra^c \la n\ra^s} 
\right) \\ 
&=\sup\limits_{\|\varphi\|_{L^2_{n,\tau}}\leq 1} 
\sum\limits_{n\in\Z}\int d\tau \sum\limits_{n_1\in\Z} \int d\tau_1 
\frac{\la \tau + n^2\ra^{-a}\la n \ra^k 
g(n_1,\tau_1)f(n-n_1,\tau-\tau_1) \bar {\varphi} (\tau,n)} {\la 
\tau_1 - n_1^3\ra^c\la n_1\ra^s 
\la \tau-\tau_1 + {(n-n_1)}^2\ra^b\la n-n_1\ra^k} \\ 
&=\sum\int\sum\int_{(n,n_1,\tau,\tau_1)\in\mathcal{R}_0} + 
\sum\int\sum\int_{(n,n_1,\tau,\tau_1)\in\mathcal{R}_1} + 
\sum\int\sum\int_{(n,n_1,\tau,\tau_1)\in\mathcal{R}_2} \\ 
&\equiv W_0 + W_1 + W_2, 
\end{split} 
\end{equation} 
whenever $\Z^2\times\R^2=\mathcal{R}_0\cup\mathcal{R}_1\cup\mathcal{R}_2$. 
 
Now, taking into account the previous calculation, we look at 
three general simple ways to reduce the problem of goods bounds on 
the expressions $W_i$ into some multiplier estimates. In the 
sequel, $\chi_{\mathcal R}$ denotes the characteristic function of 
the set ${\mathcal R}$. So, we consider the expression 
\begin{equation}\label{e.w} 
W = \sup\limits_{\|\varphi\|_{L^2_{n,\tau}}\leq 1} 
\sum\limits_{n\in\Z}\int d\tau \sum\limits_{n_1\in\Z} \int d\tau_1 
\frac{\la \tau + n^2\ra^{-a}\la n \ra^k 
g(n_1,\tau_1)f(n-n_1,\tau-\tau_1) \bar {\varphi} 
(\tau,n)\chi_{\mathcal R}} {\la \tau_1 - n_1^3\ra^b\la n_1\ra^s 
\la \tau-\tau_1 + {(n-n_1)}^2\ra^b\la n-n_1\ra^k}. 
\end{equation} 
 
The first way to bound $W$ is: integrate over $\tau_1$ and $n_1$, 
and then use the Cauchy-Schwarz and H\"older inequalities to 
obtain 
\begin{equation}\label{e.w0}\begin{split} 
|W|^2&\leq \|\varphi \|^2_{{L^2_{\tau}L^2_{n}}} \left \|\frac{\la 
n \ra^k}{\la \tau + n^2 \ra^a} \Int 
\frac{g(n_1,\tau_1)f(n-n_1,\tau-\tau_1)\chi_{{\mathcal 
R}}d\tau_1dn_1} {\la \tau_1 - n_1^3\ra^c\la n_1\ra^s\la 
\tau-\tau_1 + {(n-n_1)}^2\ra^b\la n-n_1\ra^k} 
\right \|^2_{{L^2_{\tau}L^2_{n}}}\\ 
&\leq \Int \frac{\la n \ra^{2k}}{\la \tau + n^2 \ra^{2a}} \left | 
\Int \frac{g(n_1,\tau_1)f(n-n_1,\tau-\tau_1)\chi_{{\mathcal 
R}}d\tau_1dn_1} {\la \tau_1 - n_1^3\ra^c\la n_1\ra^s\la 
\tau-\tau_1 + 
{(n-n_1)}^2\ra^b\la n-n_1\ra^k} \right |^2d\tau dn\\ 
&\leq \Int \frac{\la n \ra^{2k}}{\la \tau + n^2 \ra^{2a}} \Biggl 
(\Int \frac{\chi_{{\mathcal R}}d\tau_1dn_1} {\la \tau_1 - 
n_1^3\ra^{2c}\la n_1\ra^{2s}\la \tau-\tau_1 + 
{(n-n_1)}^2\ra^{2b}\la n-n_1\ra^{2k}}\times \\ 
&\quad \times \Int 
|g(n_1,\tau_1)|^2|f(n-n_1,\tau-\tau_1)|^2d\tau_1dn_1 
\Biggl )d\tau dn\\ 
&\leq \|f\|^2_{{L^2_{\tau}L^2_{n}}}\|g\|^2_{{L^2_{\tau_1}L^2_{n_1}}}\\ 
&\quad \times \left \| \frac{\la n \ra^{2k}}{\la \tau + n^2 
\ra^{2a}} \Int \frac{\chi_{{\mathcal R}}d\tau_1dn_1} {\la \tau_1 - 
n_1^3\ra^{2c}\la n_1\ra^{2s}\la \tau-\tau_1 + {(n-n_1)}^2\ra^{2b} 
\la n-n_1\ra^{2k}} 
\right \|_{{L^{\infty}_{\tau}L^{\infty}_{n}}}\\ 
&=\|u\|^2_{X^{k,b}}\|v\|^2_{Y^{s,c}}\\ 
&\quad \times \left \| \frac{\la n \ra^{2k}}{\la \tau + n^2 
\ra^{2a}} \Int \frac{\chi_{{\mathcal R}}d\tau_1dn_1} {\la \tau_1 - 
n_1^3\ra^{2c}\la n_1\ra^{2s} \la \tau-\tau_1 + 
{(n-n_1)}^2\ra^{2b}\la n-n_1\ra^{2k}} \right 
\|_{{L^{\infty}_{\tau}L^{\infty}_{n}}} . 
\end{split} 
\end{equation}  

The second way we can bound $W$ is: put $\widetilde f(n,\tau)= 
f(-n,-\tau)$,\; integrate over $\tau$ and $n$ first and follow the 
same steps as above to get 
\begin{equation}\label{e.w1} 
\begin{split} 
|W|^2&\leq \|g\|^2_{{L^2_{\tau_1}L^2_{n_1}}} \left \|\frac{1}{\la 
n_1\ra^s \la \tau_1 - n_1^3 \ra^c} \Int \frac{\la n \ra^k 
\widetilde 
f(n_1-n,\tau_1-\tau)\bar{\varphi}(\tau,n)\chi_{{\mathcal R}}d\tau dn}
{\la \tau + n^2\ra^a\la \tau-\tau_1 + {(n-n_1)}^2\ra^b\la n-n_1\ra^k} 
\right\|^2_{{L^2_{\tau_1}L^2_{n_1}}}\\ &\leq 
\|\widetilde f\|^2_{{L^2_{\tau_1}L^2_{n_1}}}
\|g\|^2_{{L^2_{\tau_1}L^2_{n_1}}}\\ &\quad \times 
\left\|\frac{1}{\la n_1\ra^{2s} \la \tau_1 - n_1^3 \ra^{2c}} \Int 
\frac{\la n \ra^{2k} \chi_{{\mathcal R}}d\tau dn}{\la \tau + n^2\ra^{2a}\la 
\tau-\tau_1 + {(n-n_1)}^2\ra^{2b}\la n-n_1\ra^{2k}} 
\right\|^2_{{L^{\infty}_{\tau_1}L^{\infty}_{n_1}}}\\ &= 
\|u\|^2_{X^{k,b}}\|v\|^2_{Y^{s,c}}\\ &\quad \times 
\left\|\frac{1}{\la n_1\ra^{2s} \la \tau_1 - n_1^3 \ra^{2c}} \Int 
\frac{\la n \ra^{2k} \chi_{{\mathcal R}}d\tau dn}{\la \tau + n^2\ra^{2a}\la 
\tau-\tau_1 + {(n-n_1)}^2\ra^{2b}\la n-n_1\ra^{2k}} 
\right\|_{{L^{\infty}_{\tau_1}L^{\infty}_{n_1}}}. 
\end{split} 
\end{equation}  
Note that $\widetilde f(n,\tau)=\la n \ra^k \la 
\tau - n^2\ra^b\widehat u(-n,-\tau)$\; and \; 
$\|\widetilde f\|_{{L^2_{\tau}L^2_{n}}}= \|f\|_{{L^2_{\tau}L^2_{n}}}=
\|u\|_{X^{k,b}}$.\\  Finally, the third way to estimate $W$ is: using the 
change of variables $\tau = \tau_1 -\tau_2$ and $n = n_1 -n_2$, the 
region,\;${\mathcal R}$,\;is transformed into the set $\widetilde {\mathcal R}$ 
such that 
\begin{equation*} 
\widetilde {\mathcal R} = 
\bigl \{(n_1,n_2,\tau_1,\tau_2)\in \Z^2\times\R^2;\; (n_1-n_2,n_1,
\tau_1-\tau_2,\tau_1)\in\mathcal{R} \bigl \}. 
\end{equation*}  
Then, 
$W$ can be estimated as: 
\begin{equation}\label{e.w2} 
\begin{split} |W|^2&\leq \|\widetilde{f}\|^2_{{L^2_{\tau_2}L^2_{n_2}}}\\ 
&\quad \times \left \|\frac{1}{ \la n_2\ra^k\la \tau_2 - n_2^2 \ra^b} \Int 
\frac{\la n_1 - n_2 \ra^k g(n_1,\tau_1) 
\widetilde{\bar{\varphi}}(n_2-n_1,\tau_2-\tau_1) 
\chi_{\widetilde {\mathcal R}}d\tau_1 dn_1}
{\la \tau_1-\tau_2 + (n_1-n_2)^2\ra^a\la \tau_1 - n_1^3\ra^c\la n_1\ra^s} 
\right\|^2_{{L^2_{\tau_2}L^2_{n_2}}}\\ &\leq 
\|\widetilde{f}\|^2_{{L^2_{\tau_2}L^2_{n_2}}}
\|g\|^2_{{L^2_{\tau_1}L^2_{n_1}}}\\ &\quad \times 
\left\|\frac{1}{\la n_2\ra^{2k} \la \tau_2 - n_2^2 \ra^{2b}} 
\Int \frac{\la n_1 - n_2 \ra^{2k}\chi_{\widetilde {\mathcal R}}d\tau_1 dn_1}
{\la \tau_1-\tau_2 + (n_1-n_2)^2\ra^{2a} \la \tau_1 - n_1^3\ra^{2c}
\la n_1\ra^{2s}} \right \|^2_{{L^{\infty}_{\tau_2}L^{\infty}_{n_2}}}\\ 
&= \|u\|^2_{X^{k,b}}\|v\|^2_{Y^{s,c}}\\ &\quad \times 
\left \|\frac{1}{\la n_2\ra^{2k} \la \tau_2 - n_2^2 \ra^{2b}} \Int \
\frac{\la n_1 - n_2 \ra^{2k}\chi_{\widetilde {\mathcal R}}d\tau_1 dn_1} 
{\la \tau_1-\tau_2 + (n_1-n_2)^2\ra^{2a} \la \tau_1 - n_1^3\ra^{2c}\la 
n_1\ra^{2s}} \right \|_{{L^{\infty}_{\tau_2}L^{\infty}_{n_2}}} 
\end{split} 
\end{equation}  

Next, using the equation~(\ref{e.sharp-uv-1-duality}) and the 
estimates~(\ref{e.w0}),~(\ref{e.w1}),~(\ref{e.w2}), we are going to reduce the 
desired bilinear estimate $\|uv\|_{Z^k}\lesssim 
\|u\|_{X^{k,\frac{1}{2}-}}\|v\|_{Y^{s,\frac{1}{2}}} + 
\|u\|_{X^{k,\frac{1}{2}}}\|v\|_{Y^{s,\frac{1}{2}-}}$ (whenever $s\geq 0$ and 
$k-s\leq 3/2$) into certain $L^{\infty}$ bounds for multipliers localized in 
some well-chosen regions $\mathcal{R}_i$, $i=0,1,2$ such that 
$\mathcal{R}_0\cup\mathcal{R}_1\cup\mathcal{R}_2=\Z^2\times\R^2$.  First, if 
$n_0:=n$, $n_1$, $n_2:= n_1-n$ are the frequencies of our waves, let 
$\lambda_0=\tau+n^2$, $\lambda_1:=\tau_1-n_1^3$, 
$\lambda_2:=\tau_2-n_2^2:=(\tau_1-\tau)-n_2^2$ be the modulations of our waves. 
Also, we consider $N_j = |n_j|,j=0,1,2$ variables measuring the magnitude of 
frequencies of the waves, and $L_j = |\lambda_j|,j=0,1,2$ variables measuring 
the magnitude of modulations of the waves. It is convenient to define the 
quantities $N_{max}\geq N_{med}\geq N_{min}$ to be the maximum, median and 
minimum of $N_0,N_1,N_2$, resp. Similarly, we define $L_{max}\geq L_{med}\geq 
L_{min}$.  In order to define the regions $\mathcal{R}_i$, we split 
$\Z^2\times\R^2$ into three regions ${\mathcal A}$,\;${\mathcal B}$ and 
${\mathcal C}$, 
\begin{equation*} 
\begin{split} 
&{\mathcal A}=
\bigl\{(n,n_1,\tau,\tau_1)\in \Z^2\times\R^2 ;\; N_1\leq 100 \bigl\}, \\ 
&{\mathcal B}=\bigl\{(n,n_1,\tau,\tau_1)\in \Z^2\times\R^2 ;\; 
N_1> 100\;\text{and, either}\; N_1\ll N_0\;\text{or}\;N_1\gg N_0\}, \\ 
&{\mathcal C}=\bigl\{(n,n_1,\tau,\tau_1)\in \Z^2\times\R^2 ;\; N_1> 100 
\;\text{and}\;N_1\sim N_0\}. 
\end{split} 
\end{equation*}  

Now we separate ${\mathcal C}$ into three parts 
\begin{equation*}\begin{split} 
&{\mathcal C}_0 =\bigl\{(n,n_1,\tau,\tau_1)\in {\mathcal C};\; L_0=L_{max} 
\bigl\},\\ &{\mathcal C}_1 =\bigl\{(n,n_1,\tau,\tau_1)\in {\mathcal C};\; 
L_1=L_{max}\bigl\},\\ &{\mathcal C}_2 =\bigl\{(n,n_1,\tau,\tau_1)\in 
{\mathcal C};\; L_2=L_{max}\bigl\}. 
\end{split}
\end{equation*}  

At this point, we define the sets 
${\mathcal R}_i,\;i=0,1,2$,\;as: 
\begin{equation*} 
{\mathcal R}_0= {\mathcal A}\cup {\mathcal B}\cup {\mathcal C}_0,\;\; 
{\mathcal R}_1={\mathcal C}_1,\;\; {\mathcal R}_2={\mathcal C}_2 
\end{equation*} 
and it is clear that $\Z^2\times\R^2 = {\mathcal R}_0\cup {\mathcal R}_1\cup 
{\mathcal R}_2$.  For these regions $\mathcal{R}_i$, we can show the following 
multiplier estimates  

\begin{claim}\label{c.w0} If $s\geq 0$ and $k-s\leq 3/2$, $$\left\| 
\frac{\la n \ra^{2k}}{\la \tau + n^2 \ra} \Int 
\frac{\chi_{{\mathcal R}_0}d\tau_1dn_1}{\la \tau_1 - n_1^3\ra^{1-}\la 
n_1\ra^{2s} \la \tau-\tau_1 + {(n-n_1)}^2\ra^{1-}\la n-n_1\ra^{2k}} 
\right \|_{{L^{\infty}_{\tau}L^{\infty}_{n}}}\lesssim 1.$$ 
\end{claim}  

\begin{claim}\label{c.w1} If $s\geq 0$ and $k-s\leq 3/2$, 
$$\left \|\frac{1}{\la n_1\ra^{2s} \la \tau_1 - n_1^3 \ra} \Int 
\frac{\la n \ra^{2k} \chi_{{\mathcal R}_1}d\tau dn} {\la \tau + n^2\ra\la 
\tau-\tau_1 + {(n-n_1)}^2\ra^{1-}\la n-n_1\ra^{2k}} 
\right \|_{{L^{\infty}_{\tau_1}L^{\infty}_{n_1}}}\lesssim 1.$$ 
\end{claim}  

\begin{claim}\label{c.w2} If $s\geq 0$ and $k-s\leq 3/2$, 
$$\left \|\frac{1}{\la n_2\ra^{2k} \la \tau_2 - n_2^2 \ra} \Int 
\frac{\la n_1 - n_2 \ra^{2k}\chi_{\widetilde {\mathcal R}_2}d\tau_1 dn_1} 
{\la \tau_1-\tau_2 + (n_1-n_2)^2\ra \la \tau_1 - n_1^3\ra^{1-}\la n_1\ra^{2s}} 
\right \|_{{L^{\infty}_{\tau_2}L^{\infty}_{n_2}}}\lesssim 1,$$ where 
$\widetilde {\mathcal R}_2$ is the image of $\mathcal{R}_2$ by the change of 
variables $n_2:=n_1-n$, $\tau_2:=\tau_1-\tau$. 
\end{claim}  

It is easy to show that these facts implies the desired bilinear 
estimate~(\ref{e.sharp-uv-1}). Indeed, by the 
equations~(\ref{e.w0}),~(\ref{e.w1}),~(\ref{e.w2}), we see that, for $a=1/2$ 
and well-chosen $b,c$, these claims means that, whenever $s\geq 0$ and 
$k-s\leq 3/2$, $|W_0|\lesssim 
\|u\|_{X^{k,\frac{1}{2}-}}\|v\|_{Y^{s,\frac{1}{2}-}}$, 
$|W_1|\lesssim \|u\|_{X^{k,\frac{1}{2}-}}\|v\|_{Y^{s,\frac{1}{2}}}$ and 
$|W_2|\lesssim \|u\|_{X^{k,\frac{1}{2}}}\|v\|_{Y^{s,\frac{1}{2}-}}$.  Putting 
these informations into the equation~(\ref{e.sharp-uv-1-duality}), we obtain 
the bilinear estimate~(\ref{e.sharp-uv-1}). So, it remains only to prove these 
claims.  For later use, recall the following algebraic relation: 
\begin{equation}\label{e.uv-dispersion} 
\lambda_0-\lambda_1+\lambda_2 = n_1^3-n_1^2-2 n n_1. 
\end{equation} 
\begin{proof}[Proof of claim~\ref{c.w0}] In the region $\mathcal{A}$, using 
that $N_1\leq 100$ and $\la n \ra\leq \la n_1 \ra \la n-n_1 \ra$, 
\begin{equation*} 
\begin{split} 
&\left \| \frac{\la n \ra^{2k}}{\la \tau + n^2 \ra} \sum\limits_{n_1}\int 
\frac{\chi_{{\mathcal A}}d\tau_1} {\la \tau_1 - n_1^3\ra^{1-}\la n_1\ra^{2s} 
\la \tau-\tau_1 + {(n-n_1)}^2\ra^{1-}\la n-n_1\ra^{2k}} 
\right \|_{{L^{\infty}_{\tau}L^{\infty}_{n}}} \\ &\lesssim 
\sup\limits_{n,\tau}\sum\limits_{n_1}\int \frac{d\tau_1}{\la \tau_1 - 
n_1^3\ra^{1-} \la \tau-\tau_1 + {(n-n_1)}^2\ra^{1-}} 
\end{split} 
\end{equation*} 
However, the lemma~\ref{l.calculus-1} (with $\theta=\widetilde{\theta}=1-$) 
implies 
\begin{equation*} 
\int \frac{d\tau_1}{\la \tau_1 - n_1^3\ra^{1-} \la \tau-\tau_1 + 
{(n-n_1)}^2\ra^{1-}} \leq \frac{\log (1+\la p(n_1) \ra)}{\la p(n_1) \ra^{1-}}, 
\end{equation*} 
where $p(x)$ is the polynomial $p(x):= x^3-x^2+2nx-(\tau+n^2)$. Hence, we can 
estimate: 
\begin{equation*} 
\sum\limits_{n_1}\int \frac{d\tau_1}{\la \tau_1 - n_1^3\ra^{1-} 
\la \tau-\tau_1 + {(n-n_1)}^2\ra^{1-}} \lesssim 
\sum\limits_{n_1}\frac{\log (1+\la p(n_1) \ra)}{\la p(n_1) \ra^{1-}}. 
\end{equation*} 
In particular, the lemma~\ref{l.1.1} can be applied to give 
\begin{equation}\label{e.c.w0-1} 
\left \| \frac{\la n \ra^{2k}}{\la \tau + n^2 \ra^{2a}} 
\sum\limits_{n_1}\int \frac{\chi_{{\mathcal A}}d\tau_1} 
{\la \tau_1 - n_1^3\ra^{2b}\la n_1\ra^{2s} 
\la \tau-\tau_1 + {(n-n_1)}^2\ra^{2b}\la n-n_1\ra^{2k}} 
\right \|_{{L^{\infty}_{\tau}L^{\infty}_{n}}} \lesssim 1. 
\end{equation} 
In the region $\mathcal{B}$, $N_1>100$, and either $N_1\gg N_0$ or 
$N_1\ll N_0$. In any case, it is not difficult to see that 
\begin{equation*} 
\frac{\la n \ra^{2k}}{\la n-n_1 \ra^{2k} \la n_1 \ra^{2s}}\lesssim 1. 
\end{equation*} 
In fact, this is an easy consequence of $s\geq 0$ and $N_2\gtrsim N_i$ if 
$N_i\gg N_j$, for $\{i,j\}=\{0,1\}$. So, we obtain the bound 
\begin{equation}\label{e.c.w0-2} 
\begin{split} &\left \| \frac{\la n \ra^{2k}}{\la \tau + n^2 \ra} 
\sum\limits_{n_1}\int \frac{\chi_{{\mathcal B}}d\tau_1} 
{\la \tau_1 - n_1^3\ra^{1-}\la n_1\ra^{2s} \la \tau-\tau_1 + 
{(n-n_1)}^2\ra^{1-}\la n-n_1\ra^{2k}} 
\right \|_{{L^{\infty}_{\tau}L^{\infty}_{n}}} \\ &\lesssim 
\sup\limits_{n,\tau}\sum\limits_{n_1}\int \frac{d\tau_1}{\la \tau_1 - 
n_1^3\ra^{1-} \la \tau-\tau_1 + {(n-n_1)}^2\ra^{1-}} \\ &\lesssim 
\sum\limits_{n_1}\frac{\log (1+\la p(n_1) \ra)}{\la p(n_1) \ra^{1-}} \\ 
&\lesssim 1, 
\end{split} 
\end{equation} 
where, as before, we have used the lemmas~\ref{l.calculus-1} and~\ref{l.1.1}.  

In the region $\mathcal{C}_0$, it is convenient to consider the following bound 
\begin{equation*} 
\begin{split} 
&\left \| \frac{\la n \ra^{2k}}{\la \tau + n^2 \ra} \sum\limits_{n_1}\int 
\frac{\chi_{{\mathcal C}_0}d\tau_1} {\la \tau_1 - n_1^3\ra^{1-}\la n_1\ra^{2s} 
\la \tau-\tau_1 + {(n-n_1)}^2\ra^{1-}\la n-n_1\ra^{2k}} 
\right \|_{{L^{\infty}_{\tau}L^{\infty}_{n}}} \\ 
&\lesssim \left \| \frac{1}{\la \tau + n^2 \ra} \sum\limits_{n_1}\int 
\frac{\la n_1\ra^{2k-2s}\chi_{{\mathcal C}_0}d\tau_1} 
{\la \tau_1 - n_1^3\ra^{1-}\la \tau-\tau_1 + {(n-n_1)}^2\ra^{1-}} 
\right \|_{{L^{\infty}_{\tau}L^{\infty}_{n}}}, 
\end{split} 
\end{equation*} 
which is an immediate corollary of $\la n \ra \leq \la n-n_1 \ra \la n_1 \ra$. 
Integrating with respect to $\tau_1$ and using the lemma~\ref{l.calculus-1} 
gives, as before, 
\begin{equation*} 
\begin{split} 
&\frac{1}{\la \tau + n^2 \ra} \sum\limits_{n_1}\int 
\frac{\la n_1\ra^{2k-2s}\chi_{{\mathcal C}_0}d\tau_1} 
{\la \tau_1 - n_1^3\ra^{1-}\la \tau-\tau_1 + {(n-n_1)}^2\ra^{1-}} \\ 
&\lesssim \frac{1}{\la \tau + n^2 \ra} \sum\limits_{n_1} 
\frac{\la n_1\ra^{2k-2s}\chi_{{\mathcal C}_0}\log(1+\la p(n_1) \ra)}{\la 
p(n_1) \ra^{1-}} 
\end{split} 
\end{equation*} 
Since, by the dispersion relation~(\ref{e.uv-dispersion}), 
$L_0=L_{max}\gtrsim N_1^3$ in the region $\mathcal{C}_0$, we have 
\begin{equation*} 
\begin{split} 
&\frac{1}{\la \tau + n^2 \ra} \sum\limits_{n_1} 
\frac{\la n_1\ra^{2k-2s}\chi_{{\mathcal C}_0} 
\log(1+\la p(n_1)\ra)}{\la p(n_1) \ra^{1-}} \\ 
&\lesssim \frac{L_{max}^{(2k-2s)/3}}{L_{max}}\sum\limits_{n_1} 
\frac{\log(1+\la p(n_1) \ra)}{\la p(n_1) \ra^{1-}} 
\end{split} 
\end{equation*}  
Hence, $k-s\leq 3/2$ and lemma~\ref{l.1.1} together allow us to conclude 
\begin{equation}\label{e.c.w0-3} 
\left \| \frac{\la n \ra^{2k}}{\la \tau + n^2 \ra} \sum\limits_{n_1}\int 
\frac{\chi_{{\mathcal C}_0}d\tau_1} {\la \tau_1 - n_1^3\ra^{1-}\la n_1\ra^{2s} 
\la \tau-\tau_1 + {(n-n_1)}^2\ra^{1-}\la n-n_1\ra^{2k}} 
\right \|_{{L^{\infty}_{\tau}L^{\infty}_{n}}} \lesssim 1. 
\end{equation} 
By definition of $\mathcal{R}_0$, 
the bounds~(\ref{e.c.w0-1}),~(\ref{e.c.w0-2}),~(\ref{e.c.w0-3}) concludes the 
proof of the claim~\ref{c.w0}. 
\end{proof}  

\begin{proof}[Proof of the claim~\ref{c.w1}] Using that $\la n \ra\leq \la n_1 
\ra \la n-n_1 \ra$, integrating in the variable $\tau$ and applying the 
lemma~\ref{l.calculus-1} (with $\theta=1/2$ and 
$\widetilde{\theta}=\frac{1}{2}-$), we get 
\begin{equation*} 
\begin{split} 
&\left \|\frac{1}{\la n_1\ra^{2s} \la \tau_1 - n_1^3 \ra} 
\sum\limits_{n}\int \frac{\la n \ra^{2k} \chi_{{\mathcal R}_1}d\tau } 
{\la \tau + n^2\ra\la \tau-\tau_1 + {(n-n_1)}^2\ra^{1-}\la n-n_1\ra^{2k}} 
\right \|_{{L^{\infty}_{\tau_1}L^{\infty}_{n_1}}} \\ 
&\lesssim \left \|\frac{\la n_1\ra^{2k-2s}}{\la \tau_1 - n_1^3 \ra} 
\sum\limits_{n} \frac{ \chi_{{\mathcal R}_1}\log(1+\la q(n) \ra)}
{\la q(n) \ra^{1-}} \right \|_{{L^{\infty}_{\tau_1}L^{\infty}_{n_1}}}, 
\end{split} 
\end{equation*} 
where $q(x):=\tau_1-n_1^2+2 n_1 x$. Note that in the region $\mathcal{R}_1$, 
$N_1>100$ and $N_0\sim N_1$, $|\lambda_1|\sim L_1 = L_{max}$ and, by the 
dispersion relation~(\ref{e.uv-dispersion}), $L_{max}\gtrsim N_1^3$; this 
permits us to apply the lemma~\ref{l.1.2} to conclude  
\begin{equation*} 
\frac{\la n_1\ra^{2k-2s}}{\la \tau_1 - n_1^3 \ra^{2b}} 
\sum\limits_{n} \frac{ \chi_{{\mathcal R}_1}\log(1+\la q(n) \ra)} 
{\la q(n) \ra^{1-}}\lesssim \frac{L_{max}^{(2k-2s)/3}}{L_{max}} 
\end{equation*} 
Thus, if we remember that $k-s\leq 3/2$, we get 
\begin{equation*} 
\frac{\la n_1\ra^{2k-2s}}{\la \tau_1 - n_1^3 \ra^{2b}} \sum\limits_{n} 
\frac{ \chi_{{\mathcal R}_1}\log(1+\la q(n) \ra)} 
{\la q(n) \ra^{1-}}\lesssim 1. 
\end{equation*} 
This completes the proof of the claim~\ref{c.w1}. 
\end{proof}  

\begin{proof}[Proof of the claim~\ref{c.w2}]Using that 
$\la n_1-n_2 \ra\leq \la n_1 \ra\la n_2 \ra$, integrating in the $\tau_1$ and 
applying the lemma~\ref{l.calculus-1} with $\theta=\frac{1}{2}-$ and 
$\widetilde{\theta}=1/2$,  
\begin{equation*} 
\begin{split} 
&\left \|\frac{1}{\la n_2\ra^{2k} \la \tau_2 - n_2^2 \ra^{2b}} 
\sum\limits_{n_1}\int \frac{\la n_1 - n_2 \ra^{2k}\chi_{\widetilde 
{\mathcal R}_2}d\tau_1 } {\la \tau_1-\tau_2 + (n_1-n_2)^2\ra^{2a} \la \tau_1 - 
n_1^3\ra^{2b}\la n_1\ra^{2s}} \right \|_{{L^{\infty}_{\tau_2}L^{\infty}_{n_2}}} 
\\ &\lesssim \left \|\frac{1}{\la \tau_2 - n_2^2 \ra^{2b}} \sum\limits_{n_1} 
\frac{\la n_1 \ra^{2k-2s}
\chi_{\widetilde {\mathcal R}_2}\log(1+\la r(n_1) \ra)} {r(n_1)} 
\right \|_{{L^{\infty}_{\tau_2}L^{\infty}_{n_2}}}, 
\end{split} 
\end{equation*}  
where $r(x):=x^3+x^2-2 n_2 x-(\tau_2-n_2^2)$. Note that the change of 
variables $\tau = \tau_1 -\tau_2$ and $n = n_1 -n_2$ transforms the region 
${\mathcal R}_2$ into a set $\widetilde {\mathcal R}_2$ such that 
\begin{equation*} 
\widetilde {\mathcal R}_2 \subseteq \bigl \{(n_1,n_2,\tau_1,\tau_2)\in 
\Z^2\times\R^2;\; N_1>100 \quad \text{and}\quad L_2=L_{max} \bigl \}. 
\end{equation*}  
In particular, the dispersion relation~(\ref{e.uv-dispersion}) implies that 
$|\lambda_2|\sim L_2=L_{max}\gtrsim N_1^3$ in the region  
$\widetilde{\mathcal{R}_2}$. So, an application of the lemma~\ref{l.1.1}  and 
the hypothesis $k-s\leq 3/2$ yields 
\begin{equation*} 
\frac{1}{\la \tau_2 - n_2^2 \ra^{2b}} 
\sum\limits_{n_1} \frac{\la n_1 \ra^{2k-2s}\chi_{\widetilde {\mathcal R}_2}
\log(1+\la r(n_1) \ra)} {r(n_1)} \lesssim 
\frac{L_{max}^{(2k-2s)/3}}{L_{max}}\lesssim 1. 
\end{equation*} 
This concludes the proof of the claim~\ref{c.w2}. 
\end{proof}  

It remains now only to prove the second estimate~(\ref{e.sharp-uv-2}), i.e., 
\begin{equation*} 
\left\|\frac{\la n \ra^k}{\la \tau+n^2 \ra} \widehat{uv}(n,\tau)
\right\|_{L_n^2 L_{\tau}^1} \lesssim \|u\|_{X^{k,\frac{1}{2}-}}
\cdot\|v\|_{Y^{s,\frac{1}{2}}} + \|u\|_{X^{k,\frac{1}{2}}}\cdot
\|v\|_{Y^{s,\frac{1}{2}-}} 
\end{equation*}  
We can rewrite the left-hand side as 
\begin{equation*} 
\left\|\int\limits_{n=n_1+n_2} \la n\ra^k \int\limits_{\tau=\tau_1+\tau_2} 
\frac{1}{\la\tau+n^2\ra}  \widehat{u}(n_1,\tau_1) \widehat{v}(n_2,\tau_2) 
\right\|_{L_n^2 L_{\tau}^1} 
\end{equation*} 
To begin with, we split the domain of integration into three regions. Let 
$\SL  =\SL_1\cup\SL_2\cup \SL_3$, where  $$\SL_1:=\{(n,\tau,n_2,\tau_2): 
|n_2|\leq 100\},$$  $$\SL_2:=\{(n,\tau,n_2,\tau_2): |n_2|> 100 \textrm{ and } 
|n|\ll |n_2|\},$$  $$\SL_3:=\{(n,\tau,n_2,\tau_2): |n_2|> 100 \textrm{ and } 
|n|\gg |n_2|\},$$ 
$\SM:=\{(n,\tau,n_2,\tau_2): |n_2|> 100, |n|\sim |n_2| 
\textrm{ and either } |\tau_1+n_1^2| = L_{\max} \textrm{ or } |\tau_2+n_2^3| = 
L_{\max}\}$ and $\SN:=\{(n,\tau,n_2,\tau_2): |n_2|> 100, |n|\sim |n_2| 
\textrm{ and } |\tau+n^2| = L_{\max}\}$. Clearly, $\SL$, $\SM$ and $\SN$ 
completely decomposes our domain of integrations, so that, in order to 
prove~(\ref{e.sharp-uv-2}), it suffices to get the bounds 
\begin{equation}\label{e.sharp-uv-2-L} 
\begin{split} &\left\|\int\limits_{n=n_1+n_2} \frac{\la n\ra^k}{\la n_1\ra^k  
\la n_2\ra^{s}} \int\limits_{\tau=\tau_1+\tau_2} \frac{\chi_{\SL}}{  
\la\tau+n^2\ra\la\tau_1+n_1^2\ra^{\frac{1}{2}-}  
\la\tau_2-n_2^3\ra^{\frac{1}{2}-}}\widehat{u}(n_1,\tau_1)  
\widehat{v}(n_2,\tau_2)\right\|_{L_n^2 L_{\tau}^1} \\ &\lesssim 
\|u\|_{X^{0,0}} \|v\|_{Y^{0,0}} 
\end{split} 
\end{equation} 
\begin{equation}\label{e.sharp-uv-2-M} 
\begin{split} 
&\left\|\int\limits_{n=n_1+n_2} \frac{\la n\ra^k}{\la n_1\ra^k  
\la n_2\ra^{s}} \int\limits_{\tau=\tau_1+\tau_2} \frac{\chi_{\SM}}{  
\la\tau+n^2\ra}\widehat{u}(n_1,\tau_1) \widehat{v}(n_2,\tau_2) 
\right\|_{L_n^2 L_{\tau}^1} \\ &\lesssim \|u\|_{X^{0,\frac{1}{2}-}} 
\|v\|_{Y^{0,\frac{1}{2}}} + \|u\|_{X^{0,\frac{1}{2}}} 
\|v\|_{Y^{0,\frac{1}{2}-}} 
\end{split} 
\end{equation} 
\begin{equation}\label{e.sharp-uv-2-N} 
\begin{split} 
&\left\|\int\limits_{n=n_1+n_2} \frac{\la n\ra^k}{\la n_1\ra^k\la n_2\ra^{s}} 
\int\limits_{\tau=\tau_1+\tau_2} \frac{\chi_{\SN}}{\la\tau+n^2\ra} 
\widehat{u}(n_1,\tau_1) \widehat{v}(n_2,\tau_2) \right\|_{L_n^2 L_{\tau}^1} \\ 
&\lesssim \|u\|_{X^{0,\frac{1}{2}-}} \|v\|_{Y^{0,\frac{1}{2}}} + 
\|u\|_{X^{0,\frac{1}{2}}} \|v\|_{Y^{0,\frac{1}{2}-}} 
\end{split} 
\end{equation}  

To proceed further, we need to recall the following Bourgain-Strichartz 
inequalities:  
\begin{lemma}[Bourgain~\cite{Bourgain}]\label{l.Bourgain} 
$X^{0,3/8}([0,1]), Y^{0,1/3}([0,1])\subset L^4(\T\times [0,1])$. More 
precisely, $$\|\psi(t) f\|_{L_{xt}^4}\lesssim \|f\|_{X^{0,3/8}} \quad 
\textrm{ and } \quad \|\psi(t) g\|_{L_{xt}^4}\lesssim \|g\|_{Y^{0,1/3}}.$$ 
\end{lemma}  

To prove the first bound~(\ref{e.sharp-uv-2-L}), we start with the simple 
observation that $$\frac{\la n\ra^k}{\la n_1\ra^k  \la n_2\ra^{s}} 
\lesssim 1,$$ if either $|n_2|\leq 100$, or  $|n_2|>100$ and $|n|\ll |n_2|$, 
or $|n_2|>100$ and $|n|\gg|n_2|$.  This follows from the fact that 
$\la n\ra\leq \la n_1\ra\la  n_2\ra$ and $s\geq 0$. Hence, 
\begin{equation*} 
\begin{split} 
&\left\|\int\limits_{n=n_1+n_2} \frac{\la n\ra^k}{\la n_1\ra^k \la n_2\ra^{s}} 
\int\limits_{\tau=\tau_1+\tau_2} 
\frac{\chi_{\SL}}{\la\tau+n^2\ra\la\tau_1+n_1^2\ra^{\frac{1}{2}-} 
\la\tau_2-n_2^3\ra^{\frac{1}{2}-}}\widehat{u}(n_1,\tau_1) 
\widehat{v}(n_2,\tau_2)\right\|_{L_n^2 L_{\tau}^1} \\ 
&\lesssim \left\|\int\limits_{n=n_1+n_2} \int\limits_{\tau=\tau_1+\tau_2} 
\frac{1}{\la\tau+n^2\ra \la\tau_1+n_1^2\ra^{\frac{1}{2}-} 
\la\tau_2-n_2^3\ra^{\frac{1}{2}-}}\widehat{u}(n_1,\tau_1) 
\widehat{v}(n_2,\tau_2)\right\|_{L_n^2 L_{\tau}^1}. 
\end{split} 
\end{equation*}  
Therefore, this reduces our goal to prove that 
\begin{equation*} 
\left\|\int\limits_{n=n_1+n_2} \int\limits_{\tau=\tau_1+\tau_2} 
\frac{1}{\la\tau+n^2\ra \la\tau_1+n_1^2\ra^{\frac{1}{2}-} 
\la\tau_2-n_2^3\ra^{\frac{1}{2}-}} \widehat{u}(n_1,\tau_1) 
\widehat{v}(n_2,\tau_2) \right\|_{L_n^2 L_{\tau}^1} 
\lesssim \|u\|_{X^{0,0}} \|v\|_{Y^{0,0}}. 
\end{equation*}  
This can be re-written as 
\begin{equation*} 
\left\|\frac{1}{\la\tau+n^2\ra^{5/8}\la\tau+n^2\ra^{3/8}} 
\int\limits_{n=n_1+n_2} \int\limits_{\tau=\tau_1+\tau_2} 
\widehat{u}(n_1,\tau_1) \widehat{v}(n_2,\tau_2) \right\|_{L_n^2 L_{\tau}^1} 
\lesssim \|u\|_{X^{0,\frac{1}{2}-}} \|v\|_{Y^{0,\frac{1}{2}-}}. 
\end{equation*}  
Since $2(-5/8)<-1$, the Cauchy-Schwarz inequality in $\tau$ reduces this bound 
to showing 
\begin{equation*} 
\left\|\frac{1}{\la\tau+n^2\ra^{3/8}}\int\limits_{n=n_1+n_2} 
\int\limits_{\tau=\tau_1+\tau_2} \widehat{u}(n_1,\tau_1) 
\widehat{v}(n_2,\tau_2) \right\|_{L_n^2 L_{\tau}^2} \lesssim 
\|u\|_{X^{0,\frac{1}{2}-}} \|v\|_{Y^{0,\frac{1}{2}-}}. 
\end{equation*}  
However, this bound is an easy consequence of duality, 
$L^4_{xt} L^2_{xt} L^4_{xt}$ H\"older and the Bourgain-Strichartz inequalities 
$X^{0,3/8}, Y^{0,1/3}\subset L^4$ in the lemma~\ref{l.Bourgain}.  

The second bound~(\ref{e.sharp-uv-2-M}) can be proved in an analogous fashion, 
using the dispersion relation 
\begin{equation}\label{e.uv-Dispersion} 
(\tau+n^2) - (\tau_2-n_2^3) + (\tau_1+n_1^2) = n_2^3 - n_2^2 - 2 n n_2. 
\end{equation} 
which implies that, in the region $\SM$, either $|\tau_1+n_1^2|\gtrsim |n_2|^3$ 
or $|\tau_2-n_2^3|\gtrsim |n_2|^3$. Thus, using that $k-s\leq 3/2$ and making 
the corresponding cancelation, we see that it suffices to prove that 
\begin{equation*} \left\|\frac{1}{\la\tau+n^2\ra^{5/8}\la\tau+n^2\ra^{3/8}} 
\int\limits_{n=n_1+n_2} \int\limits_{\tau=\tau_1+\tau_2} 
\widehat{u}(n_1,\tau_1) \widehat{v}(n_2,\tau_2) \right\|_{L_n^2 L_{\tau}^1} 
\lesssim \|u\|_{X^{0,0}} \|v\|_{Y^{0,\frac{1}{2}-}} 
\end{equation*} 
and 
\begin{equation*} 
\left\|\frac{1}{\la\tau+n^2\ra^{5/8}\la\tau+n^2\ra^{3/8}} \int 
\limits_{n=n_1+n_2} \int\limits_{\tau=\tau_1+\tau_2} \widehat{u}(n_1,\tau_1) 
\widehat{v}(n_2,\tau_2) \right\|_{L_n^2 L_{\tau}^1} \lesssim 
\|u\|_{X^{0,\frac{1}{2}-}} \|v\|_{Y^{0,0}}. 
\end{equation*}  
Again, we use Cauchy-Schwarz to reduce these estimates to 
\begin{equation*} 
\left\|\frac{1}{\la\tau+n^2\ra^{3/8}}\int \limits_{n=n_1+n_2} 
\int\limits_{\tau=\tau_1+\tau_2} \widehat{u}(n_1,\tau_1) 
\widehat{v}(n_2,\tau_2) \right \|_{L_n^2 L_{\tau}^2} \lesssim 
\|u\|_{X^{0,0}} \|v\|_{Y^{0,\frac{1}{2}-}} 
\end{equation*} 
and 
\begin{equation*} 
\left\|\frac{1}{\la\tau+n^2\ra^{3/8}}\int \limits_{n=n_1+n_2} 
\int\limits_{\tau=\tau_1+\tau_2} \widehat{u}(n_1,\tau_1) 
\widehat{v}(n_2,\tau_2) \right\|_{L_n^2 L_{\tau}^2} \lesssim 
\|u\|_{X^{0,\frac{1}{2}-}} \|v\|_{Y^{0,0}}, 
\end{equation*} 
which follows from duality, H\"older and Bourgain-Strichartz, as above.  

Finally, the third bound~(\ref{e.sharp-uv-2-N}) requires a subdivision into 
two cases. When $|\tau_1+n_1^2|\gtrsim |n_2|^{2-}$ (resp., 
$|\tau_2-n_2^3|\gtrsim |n_2|^{2-}$), we use $\la\tau_1+n_1^2\ra^{1/8}$ 
leaving $\la\tau_1+n_1^2\ra^{3/8}$ in the denominator and $|n_2|^{k-s-}$ in the 
numerator (resp., a similar argument with $(\tau_2-n_2^3)$ instead of 
$(\tau_1+n_1^2)$, using $\la\tau_2-n_2^3\ra^{1/6}$ and leaving 
$\la\tau_2-n_2^3\ra^{1/3}$). After another cancelation using 
$|\tau+n^2|\gtrsim |n_2|^3$, we need to prove 
\begin{equation*} 
\left\|\frac{1}{\la\tau+n^2\ra^{1/2+}}\int\limits_{n=n_1+n_2} 
\int\limits_{\tau=\tau_1+\tau_2} \widehat{u}(n_1,\tau_1) 
\widehat{v}(n_2,\tau_2) \right\|_{L_n^2 L_{\tau}^1} \lesssim 
\|u\|_{X^{0,3/8}} \|v\|_{Y^{0,\frac{1}{2}-}}, 
\end{equation*} 
and 
\begin{equation*} 
\left\|\frac{1}{\la\tau+n^2\ra^{1/2+}}\int\limits_{n=n_1+n_2} 
\int\limits_{\tau=\tau_1+\tau_2} \widehat{u}(n_1,\tau_1) 
\widehat{v}(n_2,\tau_2) \right\|_{L_n^2 L_{\tau}^1} \lesssim 
\|u\|_{X^{0,\frac{1}{2}-}} \|v\|_{Y^{0,1/3}}. 
\end{equation*}  
These bounds follow again from Cauchy-Schwarz in $\tau$, duality, H\"older and 
Bourgain-Strichartz. So it remains only the case 
$|\tau_1+n_1^2|, |\tau_2-n_2^3|\ll |n_2|^{2-}$. 
In this case, the dispersion relation says that, in the region 
$\SN$, 
$$\tau+n^2 = n_2^3-n_2^2-2n n_2 - O(|n_2|^{2-}).$$  
On the other hand, the cancelation using $|\tau+n^2|\gtrsim |n_2|^3$ and 
$k-s\leq 3/2$ reduces the proof to the bound 
\begin{equation*} 
\left\|\frac{1}{\la\tau+n^2\ra^{1/2}}\int\limits_{n=n_1+n_2} 
\int\limits_{\tau=\tau_1+\tau_2} \widehat{u}(n_1,\tau_1) 
\widehat{v}(n_2,\tau_2) \chi_{\Omega(n)}(\tau+n^2) 
\right\|_{L_n^2 L_{\tau}^1} \lesssim \|u\|_{X^{0,\frac{1}{2}-}} 
\|v\|_{Y^{0,\frac{1}{2}-}}, 
\end{equation*} 
where $\Omega(n) = \{\eta\in\R: \eta = r^3-r^2-2n r + O(|r|^{2-}), 
\text{ for some } r\in\Z, |r|\sim |n| > 100\}$. Applying Cauchy-Schwarz in 
$\tau$, we can estimate the left-hand side by 
$$\left\|\left(\int \la\tau+n^2\ra^{-1} \chi_{\Omega(n)}(\tau+n^2) 
\right)^{1/2} \left\|\int\limits_{n=n_1+n_2} \int\limits_{\tau=\tau_1+\tau_2} 
\widehat{u}(n_1,\tau_1) 
\widehat{v}(n_2,\tau_2)\right\|_{L_{\tau}^2}\right\|_{L_n^2}$$  
Therefore, the point is to show 
\begin{equation}\label{e.uv-dyadic} 
\sup_{n}\left(\int \la\tau+n^2\ra^{-1} \chi_{\Omega(n)}(\tau+n^2) 
d\tau\right)\lesssim 1 
\end{equation}  
We need the following lemma: 
\begin{lemma}\label{l.uv-dyadic}  
There exists some $\delta>0$ such that, for any fixed $n\in\Z$, $|n|\gg 1$ and 
for all $M\geq 1$ dyadic, we have 
\begin{equation*} 
|\{\mu\in\R: |\mu|\sim M, \mu = r^3-r^2-2n r + O(|r|^{2-}), \text{ for some } 
r\in\Z, |r|\sim |n|\}|\lesssim M^{1-\delta}. 
\end{equation*} 
\end{lemma}  
\begin{proof}Note that the dyadic block $\{|\mu|\sim M\}$ contains at most 
$O(M/N^2)+1$ integer numbers of the form $r^3-r^2-2n r$ with $|r|\sim |n|$, 
$r\in\Z$, where $N\sim |n|$. Indeed, this follows from the fact that the 
distance between two consecutive numbers of this form is $\sim N^2$. Thus, the 
set of $\mu$ verifying $\mu = r^3-r^2-2n r + O(|r|^{2-})$ is the union of 
$O(M/N^2)+1$ intervals of size $O(N^{2-})$. Since the relation $\mu = 
r^3-r^2-2n r + O(|r|^{2-})$ with $|\mu|\sim M$ and $|r|\sim |n|\sim N\gg 1$ 
implies that $M\sim N^3$, we get  
\begin{equation*} 
\begin{split} 
&\left|\{\mu\in\R: |\mu|\sim M, \mu = r^3-r^2-2n r + O(|r|^{2-}), 
\text{ for some } r\in\Z, |r|\sim |n|\}\right|\lesssim \\ 
& N^{2-}\cdot\frac{M}{N^2} \lesssim M^{1-}. 
\end{split} 
\end{equation*}  
This completes the proof of the lemma~\ref{l.uv-dyadic} 
\end{proof}  
Using the lemma~\ref{l.uv-dyadic}, it is not difficult to conclude the proof 
of~(\ref{e.uv-dyadic}): by changing variables, we have to estimate 
$$\sup_n \int \la \mu \ra^{-1} \chi_{\Omega(n)}(\mu) d\mu.$$ 
By decomposing the domain of integration into dyadic blocks $\{|\mu|\sim M\}$, 
the lemma~\ref{l.uv-dyadic} gives 
\begin{equation*} 
\int \la \mu \ra^{-1} \chi_{\Omega(n)}(\mu) d\mu \leq 1+ 
\sum_{M\geq 1}\int\limits_{|\mu|\sim M} \la \mu \ra^{-1} 
\chi_{\Omega(n)}(\mu) d\mu \lesssim 1+\sum\limits_{M\geq 1; \, M 
\textrm{ dyadic }} M^{-1} M^{1-\delta} \lesssim 1. 
\end{equation*} 
This proves the estimate~(\ref{e.sharp-uv-2}), thus completing the proof of 
the lemma~\ref{l.sharp-uv}. 
\end{proof}  

\subsection{Proof of proposition~\ref{p.du2}: bilinear estimates for the 
coupling term $\p_x (|u|^2)$} By the lemma~\ref{l.example-du2}, it suffices to 
prove the bilinear estimate:  
\begin{lemma}\label{l.sharp-du2}$\|\p_x(u_1 \ov{u_2})\|_{W^s}\lesssim 
\|u_1\|_{X^{k,\frac{1}{2}-}}\|u_2\|_{X^{k,\frac{1}{2}}} + 
\|u_1\|_{X^{k,\frac{1}{2}}}\|u_2\|_{X^{k,\frac{1}{2}-}}$ 
whenever $1+s\leq 4k$ and $k-s\geq -1/2$. 
\end{lemma}  
\begin{proof}From the definition of $W^s$, we have to prove that 
\begin{equation}\label{e.sharp-du2-1} 
\|\p_x(u_1 \ov{u_2})\|_{Y^{s,-1/2}}\lesssim 
\|u_1\|_{X^{k,\frac{1}{2}-}}\|u_2\|_{X^{k,\frac{1}{2}}} + 
\|u_1\|_{X^{k,\frac{1}{2}}}\|u_2\|_{X^{k,\frac{1}{2}-}} 
\end{equation} 
and 
\begin{equation}\label{e.sharp-du2-2} 
\left\|\frac{\la n\ra^s}{\la \tau-n^3 \ra}\widehat{\p_x(u_1 \ov{u_2})} 
\right\|_{L_n^2 L_{\tau}^1}\lesssim 
\|u_1\|_{X^{k,\frac{1}{2}-}}\|u_2\|_{X^{k,\frac{1}{2}}} + 
\|u_1\|_{X^{k,\frac{1}{2}}}\|u_2\|_{X^{k,\frac{1}{2}-}}. 
\end{equation}  

We begin with the proof of~(\ref{e.sharp-du2-1}). First, we reduce the 
bilinear estimate to some multiplier estimates as follows. By the definition of 
Bourgain's spaces, 
\begin{equation*} 
\begin{split} 
\|\p_x (u_1\ov {u_2})\|_{Y^{s,-a}}&= \|\la \tau-n^3 \ra^{-a} \la n \ra^s 
\widehat{\p_x (u_1 \ov{u_2})}\|_{L^2_{\tau}L^2_n} \\ &= \| n \la \tau-n^3 
\ra^{-a} \la n \ra^s \widehat{u_1}*\widehat{\ov{u_2}} 
(n,\tau)\|_{L^2_{\tau}L^2_n} 
\end{split} 
\end{equation*}  
Let 
\begin{equation*} 
f(n,\tau) = \la n \ra^k \la \tau+n^2 \ra^b \widehat{u}_1(n,\tau) \quad 
\text{ and } \quad g(n,\tau) = \la n \ra^k \la -\tau+n^2 \ra^c 
\ov{\widehat{u}_2(-n,-\tau)} 
\end{equation*}  By duality, 
\begin{equation}\label{e.sharp-du2-1-duality} 
\begin{split} 
\|\p_x (u_1\ov {u_2})\|_{Y^{s,-a}}&= \sup\limits_{\|\varphi\|_{L_{\tau}^2 
l_n^2}\leq 1} \sum\limits_{n\in\Z}\int d\tau \sum\limits_{n_1\in\Z}\int 
d\tau_1 \frac{|n| \la n \ra^s}{\la \tau-n^3 \ra^{a}}\widehat{u}_1(n-n_1,
\tau-\tau_1) \ov{\widehat{u}_2(-n_1,-\tau_1)} \cdot \ov{\varphi(n,\tau)} \\ 
&= \sup\limits_{\|\varphi\|_{L_{\tau}^2 l_n^2}\leq 1} \sum\limits_{n\in\Z}\int 
d\tau \sum\limits_{n_1\in\Z}\int d\tau_1 \frac{|n| \la n \ra^s \la \tau-n^3 
\ra^{-a} f(n-n_1,\tau-\tau_1) g(n_1,\tau_1) \ov{\varphi(n,\tau)}} {\la n-n_1 
\ra^k \la (\tau-\tau_1)+(n-n_1)^2 \ra^b \la n_1 \ra^k \la -\tau_1+n_1^2 \ra^c} 
\\ &= \sum\int\sum\int_{(n,n_1,\tau,\tau_1)\in\mathcal{V}_0} + 
\sum\int\sum\int_{(n,n_1,\tau,\tau_1)\in\mathcal{V}_1} + 
\sum\int\sum\int_{(n,n_1,\tau,\tau_1)\in\mathcal{V}_2} \\ 
&\equiv V_0+V_1 + V_2 
\end{split} 
\end{equation} 
whenever 
$\Z^2\times\R^2 = \mathcal{V}_0\cup\mathcal{V}_1\cup\mathcal{V}_2$.  

As before, we have three general ways to estimate the quantity 
\begin{equation}\label{e.v} 
V=\sup\limits_{\|\varphi\|_{L_{\tau}^2 l_n^2}\leq 1} 
\sum\limits_{n\in\Z}\int d\tau \sum\limits_{n_1\in\Z}\int d\tau_1 
\frac{|n| \la n \ra^s \la \tau-n^3 \ra^{-a} f(n-n_1,\tau-\tau_1) 
g(n_1,\tau_1) \ov{\varphi(n,\tau)}} {\la n-n_1 \ra^k \la 
(\tau-\tau_1)+(n-n_1)^2 \ra^b \la n_1 \ra^k \la -\tau_1+n_1^2 
\ra^c}\chi_{\mathcal{V}} 
\end{equation}  

Firstly, we integrate over $\tau_1$ and $n_1$ and then use Cauchy-Schwarz and 
H\"older inequalities to obtain 
\begin{equation}\label{e.v0} 
\begin{split} 
|V|^2 &\leq \|\varphi\|_{L_{n,\tau}^2}^2 \left\|\frac{|n| 
\la n \ra^s}{\la \tau-n^3 \ra^{a}}\sum\limits_{n_1}\int d\tau_1 
\frac{g(n_1,\tau_1) f(n-n_1,\tau-\tau_1)\chi_{\mathcal{V}}}{\la n-n_1 
\ra^k \la (\tau-\tau_1)+(n-n_1)^2 \ra^b \la n_1 \ra^k \la -\tau_1+n_1^2 
\ra^c} \right\|_{L_{\tau}^2 L_n^2} \\ &\leq\|u_1\|_{X^{k,b}}^2 
\|u_2\|_{X^{k,c}}^2 \times \\ &\quad \times 
\left\| \frac{|n|^2 \la n \ra^{2s}}{\la \tau - n^3 \ra^{2a}} 
\sum\limits_{n_1} \frac{1}{\la n_1 \ra^{2k} \la n-n_1\ra^{2k}} 
\int d\tau_1 \frac{\chi_{\mathcal{V}}} {\la -\tau_1 + n_1^2\ra^{2c}\la 
(\tau-\tau_1) + (n-n_1)^2\ra^{2b}}\right\|_{L_{\tau}^{\infty} L_n^{\infty}}. 
\end{split} 
\end{equation}  

Secondly, we put $\widetilde{f}(n,\tau) = f(-n,-\tau)$, integrate over $n$ and 
$\tau$ and then use the same steps above to get 
\begin{equation}\label{e.v1} 
\begin{split} 
|V|^2&\leq \|g\|_{L_{\tau_1}^2 L_{n_1}^2}^2 
\left\| \frac{1}{\la n_1 \ra^k \la -\tau_1+n_1^2 \ra^c}\sum\limits_n 
\int d\tau \frac{|n| \la n \ra^s}{\la \tau-n^3 \ra^{a}} 
\frac{\widetilde{f}(n_1-n,\tau_1-\tau)\ov{\varphi(n,\tau)}\chi_{\mathcal{V}}} 
{\la (\tau-\tau_1) + (n-n_1)^2 \ra^b \la n_1-n \ra^k}
\right\|_{L_{\tau_1}^2 L_{n_1}^2}^2 \\ &\leq \|u_1\|_{X^{k,b}}^2 
\|u_2\|_{X^{k,c}}^2 \times \\ &\quad \times \left\| 
\frac{1}{\la n_1 \ra^{2k} \la -\tau_1+n_1^2 \ra^{2c}} \sum\limits_n 
\int d\tau \frac{|n|^2 \la n \ra^{2s}\chi_{\mathcal{V}}}{\la \tau-n^3 \ra^{2a} 
\la (\tau-\tau_1) + (n-n_1)^2 \ra^{2b} \la n_1-n \ra^{2k}}
\right\|_{L_{\tau_1}^{\infty} L_{n_1}^{\infty}}. 
\end{split} 
\end{equation}  

Finally, using the change of variables $\tau_1=\tau+\tau_2$ and $n_1=n+n_2$, we 
transform $\mathcal{V}$ into the region 
\begin{equation*} 
\widetilde{\mathcal{V}}= \{(n,n_2,\tau,\tau_2): (n,n+n_2,\tau,\tau+\tau_2)\in 
\mathcal{V}\} 
\end{equation*} 
and, hence, integrating over $\tau$ and $n$, we can estimate 
\begin{equation}\label{e.v2} 
\begin{split} 
|V|^2 &\leq \|\widetilde{f}\|_{L_{\tau_2}^2 L_{n_2}^2}^2 
\left\|\frac{1}{\la n_2 \ra^k \la -\tau_2 + n_2^2 \ra^b} 
\sum\limits_{n\in\Z} \int d\tau \frac{|n|\la n \ra^s g(n+n_2,\tau+\tau_2) 
\ov{\varphi(n,\tau)}\chi_{\widetilde{\mathcal{V}}}}{\la \tau - n^3 \ra^a 
\la -(\tau+\tau_2)+(n+n_2)^2 \ra^c} \right\|_{L_{\tau_2}^2 L_{n_2}^2}^2 \\ 
&\leq \|u_1\|_{X^{k,b}}^2 \|u_2\|_{X^{k,c}}^2 \times \\ 
&\times \left\|\frac{1}{\la n_2 \ra^{2k} \la -\tau_2 + n_2^2 \ra^{2b}} 
\sum\limits_{n\in\Z} \frac{|n|^2\la n \ra^{2s}}{\la n+n_2 \ra^{2k}} 
\int d\tau \frac{\chi_{\widetilde{\mathcal{V}}}}{\la \tau - n^3 \ra^{2a} 
\la -(\tau+\tau_2)+(n+n_2)^2 \ra^{2c}} 
\right\|_{L_{\tau_2}^{\infty} L_{n_2}^{\infty}} 
\end{split} 
\end{equation}  

The next step is to use the estimates~(\ref{e.v0}),~(\ref{e.v1}) 
and~(\ref{e.v2}) for the expression~(\ref{e.v}) to reduce the bilinear estimate 
$\|\p_x(u_1 \ov{u_2})\|_{Y^{s,-1/2}}\lesssim 
\|u_1\|_{X^{k,\frac{1}{2}-}}\|u_2\|_{X^{k,\frac{1}{2}}} + 
\|u_1\|_{X^{k,\frac{1}{2}}}\|u_2\|_{X^{k,\frac{1}{2}-}}$ to $L^{\infty}$ 
bounds for certain multipliers localized in some well-chosen regions 
$\mathcal{V}_0$, $\mathcal{V}_1$ and $\mathcal{V}_2$. We consider 
$n_0:=n$, $n_1$ and $n_2:=n_1-n$ the frequencies of our waves and 
$\lambda_0:=\tau-n^3$, $\lambda_1:=-\tau_1+n_1^2$, 
$\lambda_2:=-\tau_2+n_2^2:=(\tau-\tau_1)+(n-n_1)^2$ the modulations of our 
waves; again, $L_j= |\lambda_j|$ are variables measuring the magnitude of the 
modulations, $j=0,1,2$. We define $L_{max}\geq L_{med}\geq L_{min}$ to be the 
maximum, median and minimum of $L_0,L_1,L_2$.  In order to define the regions 
$\mathcal{V}_i$, we split $\Z^2\times\R^2$ into three regions  
$\mathcal{O},\mathcal{P},\mathcal{Q}$, 
\begin{equation*} 
\begin{split} 
&\mathcal{O}=\{(n,n_1,\tau,\tau_1)\in\Z^2\times\R^2: |n|\leq 100\}, \\ 
&\mathcal{P}=\{(n,n_1,\tau,\tau_1)\in\Z^2\times\R^2: |n|\geq 100 \quad 
\text{ and } \quad |n_1|\gtrsim |n|^2\}, \\ 
&\mathcal{Q}=\{(n,n_1,\tau,\tau_1)\in\Z^2\times\R^2: |n|\geq 100 \quad 
\text{ and } \quad |n_1|\ll |n|^2\}. 
\end{split} 
\end{equation*}  

Now we separate $\mathcal{Q}$ into three parts 
\begin{equation*} 
\begin{split} 
&\mathcal{Q}_0=\{(n,n_1,\tau,\tau_1)\in\mathcal{C}: L_0=L_{max}\}, \\ 
&\mathcal{Q}_1=\{(n,n_1,\tau,\tau_1)\in\mathcal{C}: L_1=L_{max}\}, \\ 
&\mathcal{Q}_2=\{(n,n_1,\tau,\tau_1)\in\mathcal{C}: L_2=L_{max}\}. 
\end{split} 
\end{equation*}  

At this point, we put 
\begin{equation*} 
\begin{split} 
&\mathcal{V}_0=\mathcal{O}\cup\mathcal{P}\cup\mathcal{Q}_0, \\ 
&\mathcal{V}_1=\mathcal{Q}_1, \\ 
&\mathcal{V}_2=\mathcal{Q}_2. 
\end{split} 
\end{equation*}  

We have the following multiplier estimates: 
\begin{claim}\label{c.v0}If $1+s\leq 4k$ and $k-s\geq -1/2$, 
$$\left\| \frac{|n|^2 \la n \ra^{2s}}{\la \tau - n^3 \ra} \sum\limits_{n_1} 
\frac{1}{\la n_1 \ra^{2k} \la n-n_1\ra^{2k}} \int d\tau_1 
\frac{\chi_{\mathcal{V}_0}} {\la -\tau_1 + n_1^2\ra^{1-}\la (\tau-\tau_1) + 
(n-n_1)^2\ra^{1-}}\right\|_{L_{\tau}^{\infty} L_n^{\infty}} \lesssim 1.$$ 
\end{claim} 
\begin{claim}\label{c.v1}If $1+s\leq 4k$ and $k-s\geq -1/2$, 
$$\left\| \frac{1}{\la n_1 \ra^{2k} \la -\tau_1+n_1^2 \ra} \sum\limits_n 
\int d\tau \frac{|n|^2 \la n \ra^{2s}\chi_{\mathcal{V}_1}}{\la \tau-n^3 \ra 
\la (\tau-\tau_1) + (n-n_1)^2 \ra^{1-} \la n_1-n \ra^{2k}}
\right\|_{L_{\tau_1}^{\infty} L_{n_1}^{\infty}}\lesssim 1.$$ 
\end{claim} 
\begin{claim}\label{c.v2}If $1+s\leq 4k$ and $k-s\geq -1/2$, 
$$\left\|\frac{1}{\la n_2 \ra^{2k} \la -\tau_2 + n_2^2 \ra} 
\sum\limits_{n\in\Z} \frac{|n|^2\la n \ra^{2s}}{\la n+n_2 \ra^{2k}} \int 
d\tau \frac{\chi_{\widetilde{\mathcal{V}_2}}}{\la \tau - n^3 \ra 
\la -(\tau+\tau_2)+(n+n_2)^2 \ra^{1-}} 
\right\|_{L_{\tau_2}^{\infty} L_{n_2}^{\infty}}\lesssim 1,$$ 
where $\widetilde{\mathcal{V}_2}$ is the image of $\mathcal{V}_2$ under the 
change of variables $n_2:=n_1-n$ and $\tau_2:=\tau_1-\tau$. 
\end{claim}  

Again, it is easy to show that these facts implies the desired bilinear 
estimate~(\ref{e.sharp-du2-1}). Indeed, by the 
equations~(\ref{e.v0}),~(\ref{e.v1}),~(\ref{e.v2}), we see that, for $a=1/2$ 
and well-chosen $b,c$, these claims means that, whenever $1+s\leq 4k$ and 
$k-s\geq -1/2$, $|V_0|\lesssim 
\|u_1\|_{X^{k,\frac{1}{2}-}}\|u_2\|_{X^{k,\frac{1}{2}-}}$, 
$|V_1|\lesssim \|u_1\|_{X^{k,\frac{1}{2}-}}\|u_2\|_{X^{k,\frac{1}{2}}}$ and 
$|V_2|\lesssim \|u_1\|_{X^{k,\frac{1}{2}}}\|u_2\|_{X^{k,\frac{1}{2}-}}$.  
Putting these informations into the equation~(\ref{e.sharp-du2-1-duality}), 
we obtain the bilinear estimate~(\ref{e.sharp-uv-1}). Hence, we have only to 
prove these claims.  For later use, we recall that our dispersion relation is 
\begin{equation}\label{e.du2-dispersion} 
\lambda_0+\lambda_1-\lambda_2 = -n^3-n^2+2n_1n 
\end{equation}  

\begin{proof}[Proof of the claim~\ref{c.v0}] In the region $\mathcal{O}$, 
using that $|n|\leq 100$, 
\begin{equation*} 
\begin{split} 
&\sup\limits_{n,\tau} \frac{|n|^2 \la n \ra^{2s}}{\la \tau - n^3 \ra} 
\sum\limits_{n_1} \frac{1}{\la n_1 \ra^{2k} \la n-n_1\ra^{2k}} \int d\tau_1 
\frac{\chi_{\mathcal{O}}} {\la -\tau_1 + n_1^2\ra^{1-}\la (\tau-\tau_1) + 
(n-n_1)^2\ra^{1-}} \\ &\lesssim \frac{1}{\la \tau - n^3 \ra} \sum\limits_{n_1} 
\frac{1}{\la n_1 \ra^{2k} \la n-n_1\ra^{2k}}
\frac{1}{\la \lambda_1-\lambda_2 \ra^{1-}}, 
\end{split} 
\end{equation*} 
by the lemma~\ref{l.calculus-1}. By the dispersion 
relation~(\ref{e.du2-dispersion}) and the 
fact $\la x+y \ra\leq \la x \ra \la y \ra$, 
we obtain the bound 
\begin{equation}\label{e.c.v0-1} 
\begin{split} 
&\sup\limits_{n,\tau} \frac{|n|^2 \la n \ra^{2s}}{\la \tau - n^3 \ra} 
\sum\limits_{n_1} \frac{1}{\la n_1 \ra^{2k} \la n-n_1\ra^{2k}} 
\int d\tau_1 \frac{\chi_{\mathcal{O}}} {\la -\tau_1 + n_1^2\ra^{1-}
\la (\tau-\tau_1) + (n-n_1)^2\ra^{1-}} \\ 
&\lesssim \sup\limits_{n\neq 0}\sum\limits_{n_1} \frac{1}{\la n_1 \ra^{2k} 
\la n-n_1 \ra^{2k} \la -n^3-n^2+2nn_1 \ra^{1-}} \\ 
&\lesssim 1, 
\end{split} 
\end{equation} 
if $k>0$.  

In the region $\mathcal{P}$, we consider the cases $n_1=(n^2+n)/2$ and 
$|n_1-(n^2+n)/2|\geq 1$. Using that $|n|\lesssim |n_1|^{1/2}$, $4k\geq 1+s$, 
the dispersion relation~(\ref{e.du2-dispersion}) and the fact that 
$\la xy \ra\gtrsim\la x \ra \la y \ra$ whenever $|x|,|y|\geq 1$, we see that 
\begin{equation*} 
\begin{split} 
&\sup\limits_{n,\tau} \frac{|n|^2 \la n \ra^{2s}}{\la \tau - n^3 \ra} 
\sum\limits_{n_1} \frac{1}{\la n_1 \ra^{2k} \la n-n_1\ra^{2k}} \int d\tau_1 
\frac{\chi_{\mathcal{P}}} {\la -\tau_1 + n_1^2\ra^{1-}\la (\tau-\tau_1) + 
(n-n_1)^2\ra^{1-}} \\ &\lesssim C + \sup\limits_{n,\tau} |n|^2 \la n \ra^{2s} 
\sum\limits_{|n_1-(n^2+n)/2|\geq 1} \frac{1}{\la n_1 \ra^{2k} \la n-n_1\ra^{2k} 
\la n \ra^{1-} \la n_1-(n^2+n)/2 \ra^{1-}}. 
\end{split} 
\end{equation*}  
Thus, 
\begin{equation}\label{e.c.v0-2} 
\begin{split} 
&\sup\limits_{n,\tau} \frac{|n|^2 \la n \ra^{2s}}{\la \tau - n^3 \ra} 
\sum\limits_{n_1} \frac{1}{\la n_1 \ra^{2k} \la n-n_1\ra^{2k}} \int d\tau_1 
\frac{\chi_{\mathcal{P}}} {\la -\tau_1 + n_1^2\ra^{1-}\la (\tau-\tau_1) + 
(n-n_1)^2\ra^{1-}} \\ &\lesssim C + 
\sup\limits_{n,\tau} |n|^{1+} \la n \ra^{2s} 
\sum\limits_{|n_1-(n^2+n)/2|\geq 1} \frac{1}{\la n_1 \ra^{2k} 
\la n-n_1\ra^{2k} \la n_1-(n^2+n)/2 \ra^{1-}} \\ &\lesssim C + 
\sum\limits_{|n_1-(n^2+n)/2|\geq 1} \frac{1}{\la n_1 \ra^{\frac{1}{2}-} 
\la n_1-(n^2+n)/2 \ra^{1-}} \\ &\lesssim 1. 
\end{split} 
\end{equation}  

In the region $\mathcal{Q}_0$, using that $L_0\gtrsim |n|^3$ and 
$k-s\geq -1/2$, we get 
\begin{equation}\label{e.c.v0-3} 
\begin{split} 
&\sup\limits_{n,\tau} \frac{|n|^2 \la n \ra^{2s}}{\la \tau - n^3 \ra} 
\sum\limits_{n_1} \frac{1}{\la n_1 \ra^{2k} \la n-n_1\ra^{2k}} \int d\tau_1 
\frac{\chi_{\mathcal{Q}_0}} {\la -\tau_1 + n_1^2\ra^{1-}\la (\tau-\tau_1) + 
(n-n_1)^2\ra^{1-}} \\ &= \sup\limits_{n,\tau} 
\frac{|n|^2 \la n \ra^{2s}}{\la \tau - n^3 \ra} \sum\limits_{|n_1|\gtrsim |n|} 
\frac{1}{\la n_1 \ra^{2k} \la n-n_1\ra^{2k}} \int d\tau_1 
\frac{\chi_{\mathcal{Q}_0}} {\la -\tau_1 + n_1^2\ra^{1-}\la (\tau-\tau_1) + 
(n-n_1)^2\ra^{1-}} + \\ &\sup\limits_{n,\tau} 
\frac{|n|^2 \la n \ra^{2s}}{\la \tau - n^3 \ra} \sum\limits_{|n_1|\ll |n|} 
\frac{1}{\la n_1 \ra^{2k} \la n-n_1\ra^{2k}} \int d\tau_1 
\frac{\chi_{\mathcal{Q}_0}} {\la -\tau_1 + n_1^2\ra^{1-}\la (\tau-\tau_1) + 
(n-n_1)^2\ra^{1-}} \\ &\lesssim \sup\limits_{n,\tau} 
\frac{|n|^2 \la n \ra^{2s}}{\la n \ra^3 \la n \ra^{2k}} 
\sum\limits_{|n_1|\gtrsim |n|} \frac{1}{\la n-n_1\ra^{2k}} 
\int d\tau_1 \frac{\chi_{\mathcal{Q}_0}} {\la -\tau_1 + n_1^2\ra^{1-}
\la (\tau-\tau_1) + (n-n_1)^2\ra^{1-}} + \\ &\sup\limits_{n,\tau} 
\frac{|n|^2 \la n \ra^{2s}}{\la n \ra^3 \la n \ra^{2k}} 
\sum\limits_{|n_1|\ll |n|} \frac{1}{\la n_1 \ra^{2k}} \int d\tau_1 
\frac{\chi_{\mathcal{Q}_0}} {\la -\tau_1 + n_1^2\ra^{1-}\la (\tau-\tau_1) + 
(n-n_1)^2\ra^{1-}} \\ &\lesssim 1, 
\end{split} 
\end{equation} 
if $k>0$. 
\end{proof}  

\begin{proof}[Proof of the claim~\ref{c.v1}] In the region $\mathcal{Q}_1$, 
using that $L_1=L_{max}\gtrsim |n|^3$ (by the dispersion 
relation~(\ref{e.du2-dispersion}) and $|n_1|\ll |n|^2$), 
$\la n \ra \leq \la n_1 \ra \la n-n_1 \ra$ and $k-s\geq -1/2$, it is not 
difficult to see that 
\begin{equation} 
\begin{split} 
&\sup\limits_{n_1,\tau_1} \frac{1}{\la n_1 \ra^{2k} \la -\tau_1+n_1^2 \ra} 
\sum\limits_n \frac{|n|^2 \la n \ra^{2s}}{\la n_1-n \ra^{2k}} 
\int d\tau \frac{\chi_{\mathcal{Q}_1}}{\la \tau-n^3 \ra \la (\tau-\tau_1) + 
(n-n_1)^2 \ra^{1-}} \\ 
&\lesssim \sup\limits_{n_1,\tau_1} \sum\limits_{n\in\Z} 
\frac{1}{\la -\tau_1 + (n-n_1)^2 +n^3 \ra^{1-}} \\ 
&\lesssim 1. 
\end{split} 
\end{equation} 
\end{proof}  

\begin{proof}[Proof of the claim~\ref{c.v2}] In the region $\mathcal{Q}_2$, 
using that $L_2=L_{max}\gtrsim |n|^3$ (by the dispersion 
relation~(\ref{e.du2-dispersion}) and $|n_1|\ll |n|^2$), 
$\la n \ra \leq \la n_2 \ra \la n+n_2 \ra$ and $k-s\geq -1/2$, it follows that 
\begin{equation} 
\begin{split} 
&\sup\limits_{n_2,\tau_2} \frac{1}{\la n_2 \ra^{2k} \la -\tau_2 + n_2^2 \ra} 
\sum\limits_{n\in\Z} \frac{|n|^2\la n \ra^{2s}}{\la n+n_2 \ra^{2k}} 
\int d\tau \frac{\chi_{\widetilde{\mathcal{Q}}_2}}{\la \tau - n^3 \ra 
\la -(\tau+\tau_2)+(n+n_2)^2 \ra^{1-}} \\ &\lesssim \sup\limits_{n_2,\tau_2} 
\sum\limits_{n\in\Z} \frac{1}{\la \tau_2 - (n+n_2)^2 +n^3 \ra^{\theta}} \\ 
&\lesssim 1. 
\end{split} 
\end{equation} 
\end{proof}  

Once~(\ref{e.sharp-du2-1}) is proved, we start the proof of the 
estimate~(\ref{e.sharp-du2-2}), that is, 
\begin{equation*} 
\left\|\frac{\la n \ra^s}{\la \tau-n^3 \ra} 
\widehat{\p_x(u_1 \ov{u_2})} \right\|_{L_n^2 L_{\tau}^1}\lesssim 
\|u_1\|_{X^{k,\frac{1}{2}-}}\|u_2\|_{X^{k,\frac{1}{2}}} + 
\|u_1\|_{X^{k,\frac{1}{2}}}\|u_2\|_{X^{k,\frac{1}{2}-}}. 
\end{equation*}  
We can rewrite the left-hand side as 
\begin{equation*} 
\left\|\int\limits_{n=n_1+n_2} |n|\la n\ra^s \int\limits_{\tau=\tau_1+\tau_2} 
\frac{1}{\la\tau-n^3\ra}  
\widehat{u_1}(n_1,\tau_1) \ov{\widehat{u_2}(-n_2,-\tau_2)} 
\right\|_{L_n^2 L_{\tau}^1} 
\end{equation*}  
To begin with, we split the domain of integration into three regions. Let 
$\Ss =\Ss_1\cup\Ss_2$, where 
$$\Ss_1:=\{(n,\tau,n_2,\tau_2): |n|\leq 100\},$$ 
$$\Ss_2:=\{(n,\tau,n_2,\tau_2): |n|> 100 \textrm{ and } |n_2|\gtrsim |n|^2\},$$ 
$\ST:=\{(n,\tau,n_2,\tau_2): |n_2|> 100, |n_2|\ll |n|^2 \textrm{ and either } 
|\tau_1+n_1^2| = L_{\max} \textrm{ or } |-\tau_2+n_2^2| = L_{\max}\}$ 
and 
$\SU:=\{(n,\tau,n_2,\tau_2): |n_2|> 100, |n|\sim |n_2| \textrm{ and } 
|\tau-n^3| = L_{\max}\}$. Clearly, $\Ss$, $\ST$ and $\SU$ completely decomposes 
our domain of integrations, so that, in order to prove~(\ref{e.sharp-du2-2}), 
it suffices to get the bounds 
\begin{equation}\label{e.sharp-du2-2-s} 
\begin{split} 
&\left\|\int\limits_{n=n_1+n_2} \frac{|n|\la n\ra^s}{\la n_1\ra^k  
\la n_2\ra^{k}} \int\limits_{\tau=\tau_1+\tau_2} \frac{\chi_{\Ss}}{  
\la\tau+n^2\ra\la\tau_1+n_1^2\ra^{\frac{1}{2}-}  
\la-\tau_2+n_2^2\ra^{\frac{1}{2}-}}\widehat{u_1}(n_1,\tau_1)  
\ov{\widehat{u_2}(-n_2,-\tau_2)}\right\|_{L_n^2 L_{\tau}^1} \\ 
&\lesssim \|u_1\|_{X^{0,0}} \|u_2\|_{X^{0,0}} 
\end{split} 
\end{equation}  

\begin{equation}\label{e.sharp-du2-2-T} 
\begin{split} 
&\left\|\int\limits_{n=n_1+n_2} \frac{|n|\la n\ra^s}{\la n_1\ra^k  
\la n_2\ra^{k}} \int\limits_{\tau=\tau_1+\tau_2} 
\frac{\chi_{\ST}}{  \la\tau-n^3\ra}\widehat{u_1}(n_1,\tau_1) 
\ov{\widehat{u_2}(-n_2,-\tau_2)} \right\|_{L_n^2 L_{\tau}^1} \\ 
&\lesssim \|u_1\|_{X^{0,\frac{1}{2}-}} \|u_2\|_{X^{0,\frac{1}{2}}} + 
\|u_1\|_{X^{0,\frac{1}{2}}} \|u_2\|_{X^{0,\frac{1}{2}-}} 
\end{split} 
\end{equation} 

\begin{equation}\label{e.sharp-du2-2-U} 
\begin{split} 
&\left\|\int\limits_{n=n_1+n_2} \frac{|n|\la n\ra^s}{\la n_1\ra^k\la 
n_2\ra^{k}} \int\limits_{\tau=\tau_1+\tau_2} 
\frac{\chi_{\SU}}{\la\tau-n^3\ra} \widehat{u_1}(n_1,\tau_1) 
\ov{\widehat{u_2}(-n_2,-\tau_2)} 
\right\|_{L_n^2 L_{\tau}^1} \\ 
&\lesssim \|u_1\|_{X^{0,\frac{1}{2}-}} 
\|u_2\|_{X^{0,\frac{1}{2}}} + \|u_1\|_{X^{0,\frac{1}{2}}} 
\|u_2\|_{X^{0,\frac{1}{2}-}} 
\end{split} 
\end{equation}  

To prove~(\ref{e.sharp-du2-2-s}), we note that 
$$\frac{|n|\la n\ra^s}{\la n_1\ra^k  \la n_2\ra^{k}} \lesssim 1,$$ 
if either $|n|\leq 100$, or  $|n|>100$ and $|n_2|\gtrsim |n|^2$, since 
$\la n\ra\leq \la n_1\ra\la  n_2\ra$ and $1+s\leq 4k$. Hence, 
\begin{equation*} 
\begin{split} 
&\left\|\int\limits_{n=n_1+n_2} \frac{|n|\la n\ra^s}{\la n_1\ra^k 
\la n_2\ra^{k}} \int\limits_{\tau=\tau_1+\tau_2} 
\frac{\chi_{\Ss}}{\la\tau-n^3\ra\la\tau_1+n_1^2\ra^{\frac{1}{2}-} 
\la-\tau_2+n_2^2\ra^{\frac{1}{2}-}}\widehat{u_1}(n_1,\tau_1) 
\ov{\widehat{u_2}(-n_2,-\tau_2)}\right\|_{L_n^2 L_{\tau}^1} \\ 
&\lesssim \left\|\int\limits_{n=n_1+n_2} \int\limits_{\tau=\tau_1+\tau_2} 
\frac{1}{\la\tau-n^3\ra \la\tau_1+n_1^2\ra^{\frac{1}{2}-} 
\la-\tau_2+n_2^2\ra^{\frac{1}{2}-}}\widehat{u}(n_1,\tau_1) 
\widehat{v}(n_2,\tau_2)\right\|_{L_n^2 L_{\tau}^1}. 
\end{split} 
\end{equation*}  
Therefore, this reduces our goal to prove that 
\begin{equation*} \left\|\int\limits_{n=n_1+n_2} 
\int\limits_{\tau=\tau_1+\tau_2} \frac{1}{\la\tau-n^3\ra 
\la\tau_1+n_1^2\ra^{\frac{1}{2}-} \la-\tau_2+n_2^2\ra^{\frac{1}{2}-}} 
\widehat{u_1}(n_1,\tau_1) \ov{\widehat{u_2}(-n_2,-\tau_2)} 
\right\|_{L_n^2 L_{\tau}^1} \lesssim \|u_1\|_{X^{0,0}} \|u_2\|_{X^{0,0}}. 
\end{equation*}  
This can be re-written as 
\begin{equation*} 
\left\|\frac{1}{\la\tau-n^3\ra^{2/3}\la\tau-n^3\ra^{1/3}} 
\int\limits_{n=n_1+n_2} \int\limits_{\tau=\tau_1+\tau_2} 
\widehat{u_1}(n_1,\tau_1) \ov{\widehat{u_2}(-n_2,-\tau_2)} 
\right\|_{L_n^2 L_{\tau}^1} \lesssim \|u_1\|_{X^{0,\frac{1}{2}-}} 
\|u_2\|_{X^{0,\frac{1}{2}-}}. 
\end{equation*} 
Since $2(-2/3)<-1$, the Cauchy-Schwarz inequality in $\tau$ reduces 
this bound to showing 
\begin{equation*} 
\left\|\frac{1}{\la\tau-n^3\ra^{1/3}}\int\limits_{n=n_1+n_2} 
\int\limits_{\tau=\tau_1+\tau_2} 
\widehat{u_1}(n_1,\tau_1) \ov{\widehat{u_2}(n_2,\tau_2)} 
\right\|_{L_n^2 L_{\tau}^2} \lesssim \|u_1\|_{X^{0,\frac{1}{2}-}} 
\|u_2\|_{X^{0,\frac{1}{2}-}}, 
\end{equation*} 
which is an easy consequence of duality, $L^4_{xt} L^2_{xt} L^4_{xt}$ 
H\"older and the Bourgain-Strichartz inequalities $X^{0,3/8}, Y^{0,1/3}\subset 
L^4$ in the lemma~\ref{l.Bourgain}.  

The second bound~(\ref{e.sharp-du2-2-T}) can be proved in an analogous fashion, 
using the dispersion relation 
\begin{equation}\label{e.du2-Dispersion} 
(\tau-n^3) - (-\tau_2+n_2^2) + (\tau_1+n_1^2) = -n^3 + n^2 + 2n n_2. 
\end{equation} 
which implies that, in the region $\SM$, either 
$|\tau_1+n_1^2|\gtrsim |n|^3$ or $|-\tau_2+n_2^2|\gtrsim |n|^3$. Thus, using 
that $s-k\leq 1/2$ and making the corresponding cancelation, we see that it 
suffices to prove that 
\begin{equation*} 
\left\|\frac{1}{\la\tau-n^3\ra^{2/3}\la\tau-n^3\ra^{1/3}} 
\int\limits_{n=n_1+n_2} \int\limits_{\tau=\tau_1+\tau_2} 
\widehat{u_1}(n_1,\tau_1) \ov{\widehat{u_2}(-n_2,-\tau_2)} 
\right\|_{L_n^2 L_{\tau}^1} \lesssim \|u_1\|_{X^{0,0}} 
\|u_2\|_{X^{0,\frac{1}{2}-}} 
\end{equation*} 
and 
\begin{equation*} 
\left\|\frac{1}{\la\tau-n^3\ra^{2/3}\la\tau-n^3\ra^{1/3}} 
\int\limits_{n=n_1+n_2} \int\limits_{\tau=\tau_1+\tau_2} 
\widehat{u_1}(n_1,\tau_1) \ov{\widehat{u_2}(-n_2,-\tau_2)} 
\right\|_{L_n^2 L_{\tau}^1} \lesssim \|u_1\|_{X^{0,\frac{1}{2}-}} 
\|u_2\|_{X^{0,0}}. 
\end{equation*}  
Again, we use Cauchy-Schwarz to reduce these estimates to 
\begin{equation*} 
\left\|\frac{1}{\la\tau-n^3\ra^{1/3}}\int\limits_{n=n_1+n_2} 
\int\limits_{\tau=\tau_1+\tau_2} 
\widehat{u_1}(n_1,\tau_1) \ov{\widehat{u_2}(-n_2,-\tau_2)} 
\right\|_{L_n^2 L_{\tau}^2} \lesssim \|u_1\|_{X^{0,0}} 
\|u_2\|_{X^{0,\frac{1}{2}-}} 
\end{equation*} 
and 
\begin{equation*} 
\left\|\frac{1}{\la\tau-n^3\ra^{1/3}}\int\limits_{n=n_1+n_2} 
\int\limits_{\tau=\tau_1+\tau_2} \widehat{u_1}(n_1,\tau_1) 
\ov{\widehat{u_2}(-n_2,-\tau_2)} \right\|_{L_n^2 L_{\tau}^2} \lesssim 
\|u_1\|_{X^{0,\frac{1}{2}-}} \|u_2\|_{X^{0,0}}, 
\end{equation*} 
which follows from duality, H\"older and Bourgain-Strichartz, as above.  

Finally, the third bound~(\ref{e.sharp-du2-2-U}) requires a subdivision into 
two cases. When $|\tau_1+n_1^2|\gtrsim |n|^{1-}$ 
(resp., $|-\tau_2+n_2^2|\gtrsim |n|^{1-}$), we use $\la\tau_1+n_1^2\ra^{1/8}$ 
leaving $\la\tau_1+n_1^2\ra^{3/8}$ in the denominator and 
$|n|^{1+s-k-}$ in the numerator (resp., the same argument with 
$(-\tau_2+n_2^2)$ instead of $(\tau_1+n_1^2)$). After another cancelation using 
$|\tau-n^3|\gtrsim |n|^3$, we need to prove 
\begin{equation*} 
\left\|\frac{1}{\la\tau-n^3\ra^{1/2+}}\int\limits_{n=n_1+n_2} 
\int\limits_{\tau=\tau_1+\tau_2} \widehat{u_1}(n_1,\tau_1) 
\ov{\widehat{u_2}(-n_2,-\tau_2)} \right\|_{L_n^2 L_{\tau}^1} 
\lesssim \|u_1\|_{X^{0,3/8}} \|u_2\|_{X^{0,\frac{1}{2}-}}, 
\end{equation*} 
and 
\begin{equation*} 
\left\|\frac{1}{\la\tau-n^3\ra^{1/2+}}\int\limits_{n=n_1+n_2} 
\int\limits_{\tau=\tau_1+\tau_2} \widehat{u_1}(n_1,\tau_1) 
\ov{\widehat{u_2}(-n_2,-\tau_2)} \right\|_{L_n^2 L_{\tau}^1} 
\lesssim \|u_1\|_{X^{0,\frac{1}{2}-}} \|u_2\|_{X^{0,3/8}}. 
\end{equation*}  
These bounds follow again from Cauchy-Schwarz in $\tau$, duality, H\"older and 
Bourgain-Strichartz. So it remains only the case 
$|\tau_1+n_1^2|, |\tau_2-n_2^3|\ll |n|^{1-}$. In this case, the dispersion 
relation says that, in the region $\SN$, 
$$\tau-n^3 = -n^3+n^2+2n n_2 - O(|n|^{1-}).$$  
On the other hand, the cancelation using $|\tau-n^3|\gtrsim |n|^3$ and 
$s-k\leq 1/2$ reduces the proof to the bound 
\begin{equation*} 
\left\|\frac{1}{\la\tau-n^3\ra^{1/2}}\int\limits_{n=n_1+n_2} 
\int\limits_{\tau=\tau_1+\tau_2} \widehat{u_1}(n_1,\tau_1) 
\ov{\widehat{u_2}(-n_2,-\tau_2)} \chi_{\widetilde{\Omega}(n)}(\tau-n^3) 
\right\|_{L_n^2 L_{\tau}^1} \lesssim \|u_1\|_{X^{0,\frac{1}{2}-}} 
\|u_2\|_{X^{0,\frac{1}{2}-}}, 
\end{equation*} 
where, 
$\widetilde{\Omega}(n) = \{\eta\in\R: \eta = n^3-n^2-2n r + O(|n|^{1-}), 
\text{ for some } r\in\Z, |r|\ll |n|^2\}$ if $|n| > 100$ and 
$\widetilde{\Omega}(n) = \emptyset$ otherwise. Applying Cauchy-Schwarz in 
$\tau$, we can estimate the left-hand side by 
$$\left\|\left(\int \la\tau-n^3\ra^{-1} \chi_{\widetilde{\Omega}(n)}(\tau-n^3) 
\right)^{1/2} 
\left\|\int\limits_{n=n_1+n_2} \int\limits_{\tau=\tau_1+\tau_2} 
\widehat{u_1}(n_1,\tau_1) 
\ov{\widehat{u_2}(-n_2,-\tau_2)}\right\|_{L_{\tau}^2}\right\|_{L_n^2}$$  
Therefore, the point is to show 
\begin{equation}\label{e.du2-dyadic} 
\sup_{n}\left(\int \la\tau-n^3\ra^{-1} 
\chi_{\widetilde{\Omega}(n)}(\tau-n^3) d\tau\right)\lesssim 1 
\end{equation} 
We need the following lemma: 
\begin{lemma}\label{l.du2-dyadic} 
There exists some $\delta>0$ such that, for any fixed $n\in\Z$, $|n|\gg 1$ and 
for all $M\geq 1$ dyadic, we have 
\begin{equation*} 
|\{\mu\in\R: |\mu|\sim M, \mu = n^3-n^2-2n r + O(|n|^{1-}), 
\text{ for some } r\in\Z, |r|\ll |n|^2\}|\lesssim M^{1-\delta}. 
\end{equation*} 
\end{lemma}  
\begin{proof}Note that the dyadic block $\{|\mu|\sim M\}$ contains at most 
$O(M/N)+1$ integer numbers of the form $n^3-n^2-2n r$ with $r\in\Z$, where 
$N\sim |n|$. Indeed, this follows from the fact that the distance between two 
consecutive numbers of this form is $\sim N$. Thus, the set of $\mu$ verifying 
$\mu = r^3-r^2-2n r + O(|r|^{2-})$ is the union of $O(M/N)+1$ intervals of size 
$O(N^{2-})$. Since the relation $\mu = n^3-n^2-2n r + O(|n|^{1-})$ with 
$|\mu|\sim M$ and $|r|\ll |n|^2\sim N^2\gg 1$ implies that $M\sim N^3$, we get 
\begin{equation*} 
\begin{split} 
&\left|\{\mu\in\R: |\mu|\sim M, \mu = 
n^3-n^2-2n r + O(|n|^{1-}), \text{ for some } r\in\Z, |r|\ll 
|n|^2\}\right|\lesssim \\ & N^{1-}\cdot\frac{M}{N} \lesssim M^{1-}. 
\end{split} 
\end{equation*}  
This completes the proof of the lemma~\ref{l.du2-dyadic} 
\end{proof}  
It is now easy to conclude the proof of~(\ref{e.uv-dyadic}): by changing 
variables, we have to estimate 
$$\sup_n \int \la \mu \ra^{-1} \chi_{\widetilde{\Omega}(n)}(\mu) d\mu.$$ 
By decomposing the domain of integration into dyadic blocks 
$\{|\mu|\sim M\}$, the lemma~\ref{l.du2-dyadic} gives 
\begin{equation*} 
\begin{split} 
\int \la \mu \ra^{-1} \chi_{\widetilde{\Omega}(n)}(\mu) d\mu \leq 1+ 
\sum_{M\geq 1}\int\limits_{|\mu|\sim M} \la \mu \ra^{-1} 
\chi_{\widetilde{\Omega}(n)(\mu)} d\mu \\ 
\lesssim 1+\sum\limits_{M\geq 1; \, M \textrm{ dyadic }} 
M^{-1} M^{1-\delta} \lesssim 1. 
\end{split} 
\end{equation*}  
This proves the estimate~(\ref{e.sharp-du2-2}), thus completing the proof of 
the lemma~\ref{l.sharp-du2}. 
\end{proof}  

\section{Local well-posedness for rough initial data}\label{s.thmA}  
This section contains the proof of the theorem~\ref{t.A} concerning the local 
well-posedness of the NLS-KdV. First of all, we observe that the 
NLS-KdV~(\ref{e.nls-kdv}) is equivalent to the integral equation 
\begin{equation*} 
u(t) = U(t)u_0 - i\int_0^t U(t-t')\{\alpha u(t')v(t') 
+\beta |u|^2 u(t') \} dt', 
\end{equation*}  
\begin{equation*} 
v(t) = V(t)v_0 + \int_0^t V(t-t')\{\gamma\p_x(|u|^2)(t') - 
\frac{1}{2}\p_x(v^2)(t')\} dt'. 
\end{equation*}  
Since we are seeking for local-in-time solutions for~(\ref{e.nls-kdv}), it 
suffices to find a fixed point $u$ for the map 
$\Phi = (\Phi_1,\Phi_2):\widetilde{X}^k([0,T])\times\widetilde{Y}^s([0,T])\to 
\widetilde{X}^k([0,T])\times\widetilde{Y}^s([0,T])$, 
\begin{equation*} 
\Phi_1(u,v)=\psi_1(t)U(t)u_0 -i\psi_T(t)\int_0^t 
U(t-t')\{ \alpha u(t')v(t') + \beta |u|^2 u(t')\} dt', 
\end{equation*}  
\begin{equation*} 
\Phi_2(u,v)=\psi_1(t)V(t)v_0 + \psi_T(t)\int_0^t 
V(t-t')\{ \gamma\p_x(|u|^2)(t') - \frac{1}{2}\p_x(v^2)(t')\} dt'. 
\end{equation*}  
From now on, our efforts are to show that $\Phi$ is a contraction of 
(a large ball of) the space 
$\widetilde{X}^k([0,T])\times\widetilde{Y}^s([0,T])$ for sufficiently small 
$T>0$. To accomplish this goal, we need the following well-known linear and 
multilineal estimates related to the cubic NLS and the KdV equations:  
\begin{lemma}[Linear estimates]\label{l.linear}It holds 
\begin{itemize} 
\item $\|\psi_1(t) U(t)u_0\|_{X^k}\lesssim \|u_0\|_{H^k}$ and 
$\left\|\psi_T(t)\int_0^t U(t-t') F(t') dt'\right\|_{X^k} 
\lesssim \|F\|_{Z^k}$; 
\item $\|\psi_1(t) V(t)v_0\|_{Y^s}\lesssim \|v_0\|_{H^s}$ and 
$\left\|\psi_T(t)\int_0^t V(t-t') G(t') dt'\right\|_{Y^s}\lesssim \|G\|_{W^s}$. 
\end{itemize} 
\end{lemma}  

\begin{lemma}[Trilinear estimate for the cubic term $|u|^2u$]\label{l.nls} 
For $k\geq 0$, we have 
$$ \|\psi(t) u v \ov{w}\|_{Z^k}\lesssim 
\|u\|_{X^{k,\frac{3}{8}}} \|v\|_{X^{k,\frac{3}{8}}} \|w\|_{X^{k,\frac{3}{8}}}. 
$$ 
\end{lemma}  

\begin{lemma}[Bilinear estimate for $\p_x(v^2)$]\label{l.kdv} For $s\geq -1/2$, 
we have 
$$ \|\psi(t)\p_x(v_1 v_2)\|_{W^s}\lesssim \|v_1\|_{Y^{s,\frac{1}{2}}} 
\|v_2\|_{Y^{s,\frac{1}{2}-}} + \|v_1\|_{Y^{s,\frac{1}{2}-}} 
\|v_2\|_{Y^{s,\frac{1}{2}}}, $$ 
if $v_1=v_1(x,t)$ and $v_2=v_2(x,t)$ are $x$-periodic functions having zero 
$x$-mean for all $t$ (i.e., $\int_{\T} v_j(x,t) dx =0$ for all $t$ and $j=1,2$). 
\end{lemma}    

\begin{remark}The zero-mean assumption in the lemma~\ref{l.kdv} above is 
crucial for some of the analysis of the multiplier associated to this bilinear 
estimate. However, in the proof of our local well-posedness result, this 
hypothesis is not restrictive by a standard argument based on the conservation 
of the mean of $v$ under the flow~(\ref{e.nls-kdv}). See the 
remark~\ref{r.zero-mean} below. 
\end{remark}  

We present the proofs of these lemmas in the appendix of this paper because 
some of these estimates are not stated as above in the literature, although 
they are contained in the works~\cite{Bourgain} and~\cite{CKSTT} for instance. 
See the section~\ref{s.appendix} below for more details.  Returning to the 
proof of the theorem~\ref{t.A}, in order to apply the lemma~\ref{l.kdv}, we 
make the following observation:  

\begin{remark}\label{r.zero-mean}The spatial mean $\int_{\T} v(t,x) dx$ is 
preserved during the evolution~(\ref{e.nls-kdv}). 
Thus, we can assume that the initial data $v_0$ has zero-mean, since otherwise 
we make the change $w= v-\int_{\T}v_0 dx$ at the expense of two harmless linear 
terms (namely, $u\int_{\T}v_0 dx$ and $\p_x v \int_{\T}v_0$). 
\end{remark}  

After this reduction, we are ready to finish the proof of the 
theorem~\ref{t.A}. Accordingly with the linear estimates 
(lemma~\ref{l.linear}), trilinear estimate for the cubic term $|u|^2 u$ 
(lemma~\ref{l.nls}), bilinear estimate for $\p_x(v^2)$ (lemma~\ref{l.kdv}) and 
the bilinear estimates for the coupling terms (propositions~\ref{p.uv} 
and~\ref{p.du2}), we obtain  
\begin{equation*} 
\begin{split} 
\|\Phi_1(u,v)\|_{\widetilde{X}^k([0,T])}&\leq C_0\|u_0\|_{H^k} + 
C_1\{\|uv\|_{Z^k}+\|u\|_{X^{k,\frac{3}{8}}([0,T])}^3\} \\ 
&\leq C_0\|u_0\|_{H^k} + 
C_1 \|u\|_{X^{k,\frac{1}{2}-}([0,T])}\|v\|_{Y^{s,\frac{1}{2}}([0,T])} + \\ 
&+ C_1\|u\|_{X^{k,\frac{1}{2}}([0,T])}\|v\|_{Y^{s,\frac{1}{2}-}([0,T])} + 
C_1 \|u\|_{X^{k,\frac{3}{8}}([0,T])}^3 
\end{split} 
\end{equation*} 
and 
\begin{equation*} 
\begin{split} 
\|\Phi_2(u,v)\|_{\widetilde{Y}^s([0,T])}&\leq C_0\|v_0\|_{H^s} + 
C_1\{\|\p_x(v^2)\|_{W^s}+\|\p_x(|u|^2)\|_{W^s}\} \\ 
&\leq C_0\|v_0\|_{H^k} + 
C_1\{ \|v\|_{Y^{s,\frac{1}{2}}}\|v\|_{Y^{s,\frac{1}{2}-}([0,T])} + 
\|u\|_{X^{k,\frac{1}{2}}}\|u\|_{X^{k,\frac{1}{2}-}([0,T])} \}, 
\end{split} 
\end{equation*} 
if $s\geq 0$, $-1/2\leq k-s\leq 3/2$ and $1+s\leq 4k$.  At this point we invoke 
the following elementary lemma concerning the stability of Bourgain's spaces 
with respect to time localization:
 
\begin{lemma}Let $X_{\tau=h(\xi)}^{s,b}:=\{f: \la\tau-h(\xi)\ra^b\la\xi\ra^s 
|\widehat{f}(\tau,\xi)|\in L^2\}$. Then, 
$$\|\psi(t)f\|_{X_{\tau=h(\xi)}^{s,b}}\lesssim_{\psi,b} 
\|f\|_{X_{\tau=h(\xi)}^{s,b}}$$ 
for any $s,b\in\R$ and, furthermore, if $-1/2<b'\leq b <1/2$, then for any 
$0<T<1$ we have 
$$\|\psi_T(t)f\|_{X_{\tau=h(\xi)}^{s,b'}}\lesssim_{\psi,b',b} T^{b-b'} 
\|f\|_{X_{\tau=h(\xi)}^{s,b}},$$ 
\end{lemma}  

\begin{proof}First of all, note that $\la\tau-\tau_0-h(\xi)\ra^{b}\lesssim_b 
\la\tau_0\ra^{|b|}\la\tau-h(\xi)\ra^{b}$, 
from which we obtain 
$$\|e^{it\tau_0}f\|_{X_{\tau=h(\xi)}^{s,b}}\lesssim_b \la\tau_0\ra^{|b|} 
\|f\|_{X_{\tau=h(\xi)}^{s,b}}.$$  
Using that $\psi(t)=\int\widehat{\psi}(\tau_0) e^{it\tau_0}d\tau_0$, we 
conclude 
$$\|\psi(t)f\|_{X_{\tau=h(\xi)}^{s,b}}\lesssim_b 
\left(\int|\widehat{\psi}(\tau_0)| \la\tau_0\ra^{|b|}\right) 
\|f\|_{X_{\tau=h(\xi)}^{s,b}}.$$  
Since $\psi$ is smooth with compact support, the first estimate follows.  

Next we prove the second estimate. By conjugation we may assume $s=0$ and, 
by composition it suffices to treat the cases $0\leq b'\leq b$ or 
$\leq b'\leq b\leq 0$. By duality, we may take $0\leq b'\leq b$. 
Finally, by interpolation with the trivial case $b'=b$, we may consider 
$b'=0$. This reduces matters to show that 
$$\|\psi_T(t)f\|_{L^2}\lesssim_{\psi,b} T^b\|f\|_{X_{\tau=h(\xi)}^{0,b}}$$ 
for $0<b<1/2$. Partitioning the frequency spaces into the cases 
$\la\tau-h(\xi)\ra\geq 1/T$ and $\la\tau-h(\xi)\leq 1/T$, we see that in the 
former case we'll have 
$$\|f\|_{X_{\tau=h(\xi)}^{0,0}}\leq T^b\|f\|_{X_{\tau=h(\xi)}^{0,b}}$$ 
and the desired estimate follows because the multiplication by $\psi$ is a 
bounded operation in Bourgain's spaces. In the latter case, by Plancherel and 
Cauchy-Schwarz 
\begin{equation*} 
\begin{split} 
\|f(t)\|_{L_x^2}&\lesssim \|\widehat{f(t)}(\xi)\|_{L_{\xi}^2} \lesssim 
\left\|\int_{\la\tau-h(\xi)\ra\leq 1/T}|\widehat{f}(\tau,\xi)|d\tau) 
\right\|_{L_{\xi}^2} \\ &\lesssim_b T^{b-1/2} 
\left\|\int\la\tau-h(\xi)\ra^{2b} |\widehat{f}(\tau,\xi)|^2 
d\tau)^{1/2}\right\|_{L_{\xi}^2} = T^{b-1/2}\|f\|_{X_{\tau=h(\xi)}^{s,b}}. 
\end{split} 
\end{equation*}  
Integrating this against $\psi_T$ concludes the proof of the lemma. 
\end{proof} 

Now, a direct application of this lemma yields 
\begin{equation*} 
\|\Phi_1(u,v)\|_{\widetilde{X}^k([0,T])}\leq C_0\|u_0\|_{H^k} + 
C_1 T^{0+}\{\|u\|_{X^{k,\frac{1}{2}}([0,T])}\|v\|_{Y^{s,\frac{1}{2}}([0,T])}  
+ \|u\|_{X^{k,\frac{1}{2}}([0,T])}^3\}. 
\end{equation*} 
and 
\begin{equation*} 
\begin{split} 
\|\Phi_2(u,v)\|_{\widetilde{Y}^s([0,T])}&\leq C_0\|v_0\|_{H^s} + 
C_1T^{0+}\{\|v\|_{Y^{s,\frac{1}{2}}([0,T])}^2 + 
\|u\|_{X^{k,\frac{1}{2}}([0,T])}^2 \}, 
\end{split} 
\end{equation*} 
if $s\geq 0$, $-1/2\leq k-s\leq 3/2$ and $1+s\leq 4k$. Hence, if $T>0$ is 
sufficiently small (depending on $\|u_0\|_{H^k}$ and $\|v_0\|_{H^s}$), we see 
that for every sufficiently large $R>0$, $\Phi$ sends the ball of radius $R$ 
of the space $\widetilde{X}^k([0,T])\times\widetilde{Y}^s([0,T])$ into itself.  
Similarly, we have that 
\begin{equation*} 
\|\Phi_1(u,v)-\Phi_1(\widetilde{u},\widetilde{v})\|_{\widetilde{X}^k} 
\lesssim T^{0+}\{\|u\|_{X^{k,\frac{1}{2}}}+\|u\|_{X^{k,\frac{1}{2}}}^2+ 
\|v\|_{Y^{s,\frac{1}{2}}}\}\{\|u-\widetilde{u}\|_{X^{k,\frac{1}{2}}}+ 
\|v-\widetilde{v}\|_{Y^{s,\frac{1}{2}}}\} 
\end{equation*} 
and 
\begin{equation*} 
\|\Phi_2(u,v)-\Phi_2(\widetilde{u},\widetilde{v})\|_{\widetilde{Y}^s([0,T])} 
\lesssim T^{0+}\{\|u\|_{X^{k,1/2}}+\|v\|_{Y^{s,1/2}}\} 
\{\|u-\widetilde{u}\|_{X^{k,1/2}}+\|v-\widetilde{v}\|_{Y^{s,1/2}}\}, 
\end{equation*} 
if $s\geq 0$, $-1/2\leq k-s\leq 3/2$ and $1+s\leq 4k$. So, up to taking $T>0$ 
smaller, we get that $\Phi$ is a contraction. This concludes the proof of the 
theorem~\ref{t.A}. 

\section{Global well-posedness in the energy space $H^1\times H^1$}
\label{s.thmB}  
This section is devoted to the proof of the theorem~\ref{t.B}. First of all, we 
recall the following conserved functionals for the NLS-KdV system 
\begin{lemma}\label{l.global}The evolution~(\ref{e.nls-kdv}) preserves the 
quantities 
\begin{itemize} 
\item $M(t):=\int_{\T} |u(t)|^2 dx$, 
\item $Q(t):=\int_{\T}\left\{\alpha v(t)^2 + 2\gamma \Im(u(t)\ov{\p_x u(t)}) 
dx\right\}$ and 
\item $E(t):=\int_{\T}\left\{\alpha\gamma v(t)|u(t)|^2 -\frac{\alpha}{6}v(t)^3 
+ \frac{\beta\gamma}{2}|u(t)|^4 + \frac{\alpha}{2}|\p_x v(t)|^2 + 
\gamma |\p_x u(t)|^2\right\} dx$. 
\end{itemize}  
In other words, $M(t)=M(0)$, $Q(t)=Q(0)$ and $E(t)=E(0)$. 
\end{lemma} 
In order to do not interrupt the proof of the global well-posedness result, we 
postpone the proof of this lemma to the appendix.  

Let $\alpha \gamma >0$ and $t >0$. From the previous lemma, we have that 
$\|u(t)\|_{L^2}=\|u_0\|_{L^2}$, and 
\begin{equation*} 
\|v(t)\|_{L^2}^2 \leq \frac{1}{|\alpha|}\left\{ |\mathcal{Q}_0| + 
2|\gamma|\;\|u_0\|_{L^2}\|\partial_x u(t)\|_{L^2}\right\}. 
\end{equation*}  
Put $\mu=\min \left \{|\gamma|,\tfrac{|\alpha|}{2}\right \}$. Then, using 
again the previous lemma, Gagliardo-Nirenberg and Young inequalities, we 
deduce  
\begin{equation*} 
\begin{split} 
\|\partial_x u(t)\|_{L^2}^2&+ \|\partial_x v(t)\|_{L^2}^2 \leq 
\frac{1}{\mu}\left( |\gamma|\|\partial_x u(t)\|_{L^2}^2 + 
|\alpha|\|\partial_x v(t)\|_{L^2}^2\right )\\ 
&\leq C\left( |E(0)|+ \|v(t)\|_{L^2}\|u(t)\|_{L^4}^2 + 
\|v(t)\|_{L^3}^3+ \|u(t)\|_{L^4}^4 \right)\\ 
&\leq C\left(|E(0)|+ \|v(t)\|_{L^2}^2 + 
\|v(t)\|_{L^3}^3+ \|u(t)\|_{L^4}^4 \right)\\ 
&\leq C\left(|E(0)|+ |Q(0)| + \|u_0\|_{L^2}\|\partial_x u(t)\|_{L^2} + 
\|v(t)\|_{L^3}^3+ \|u(t)\|_{L^4}^4 \right)\\ 
&\leq  C \left\{|E(0)|+ |Q(0)| + |Q(0)|^{\frac{5}{3}} + M(0)^5 + M(0)^3 + M(0) 
\right\} + \\ 
&+ \frac{1}{2}\left\{ \|\partial_x u(t)\|_{L^2}^2+ 
\|\partial_x v(t)\|_{L^2}^2\right\}. 
\end{split}
\end{equation*}  
Hence 
\begin{equation}
\begin{split}\label{e.global-1} 
\|\partial_x u(t)\|_{L^2}^2+\|\partial_x v(t)\|_{L^2}^2 
&\leq C\left\{|E(0)|+ |Q(0)| + |Q(0)|^{\frac{5}{3}} + M(0)^5 + M(0)^3 + M(0) 
\right\}. 
\end{split}
\end{equation}  
We can estimate the right hand of (\ref{e.global-1}) using the conservation 
laws in the lemma~\ref{l.global} and Sobolev's lemma to get 
\begin{equation}\label{e.global-2} 
\|u(t)\|_{H^1}^2+\|v(t)\|_{H^1}^2 \leq \Psi(\|u_0\|_{H^1}, \|v_0\|_{H^1}), 
\end{equation} 
where $\Psi$ is a function depending only on 
$\|u_0\|_{H^1}$\;and\;$\|v_0\|_{H^1}$.  We observe that the constans depend 
only on the parameters 
$\alpha, \beta$ and $\gamma$.  
Since the right hand in (\ref{e.global-2}) only depends of 
$\|u_0\|_{H^1}$\;and\;$\|v_0\|_{H^1}$, we can repeat the argument of local 
existence of solution at time $T$ arriving to a solution for any positive time. 
This completes the proof of the theorem~\ref{t.B}. 

\section{Final Remarks}\label{s.remarks} 
We conclude this paper with some comments and questions related to our results 
in theorems~\ref{t.A},~\ref{t.B}.  

Concerning the local well-posedness result in theorem~\ref{t.A}, the gap 
between our endpoint $H^{1/4}\times L^2$ and the ``natural'' 
$L^2\times H^{-1/2}$ endpoint\footnote{As we said before, from the sharp 
well-posedness theory for the NLS and the KdV equations, the well-posedness 
endpoint for the periodic NLS equation is $L^2$ and for the periodic KdV is 
$H^{-1/2}$.} suggests the ill-posedness question:  
\begin{question}Is the periodic NLS-KdV system~(\ref{e.nls-kdv}) ill-posed for 
initial data $(u_0,v_0)\in H^k\times H^s$ with $0\leq k<1/4$, $1+s\leq 4k$ and 
$-1/2\leq k-s\leq 3/2$? 
\end{question}  

On the other hand, one should be able to improve the global well-posedness 
result in theorem~\ref{t.B} using the \emph{I-method} of Colliander, Keel, 
Staffilani, Takaoka and Tao~\cite{CKSTT}. In the continuous case, the global 
well-posedness result in the energy space of Corcho and Linares~\cite{Corcho} 
was refined by Pecher~\cite{Pecher} via the I-method. This motivates the 
following question in the periodic context: 
\begin{question}Is the periodic NLS-KdV system~(\ref{e.nls-kdv}) globally 
well-posed for initial data $(u_0,v_0)\in H^{1-}\times H^{1-}$? 
\end{question}  
We plan to address this issue in a forthcoming paper by using our bilinear 
estimates for the coupling terms $uv$ and $\p_x(|u|^2)$ and the I-method.   

\section*{Acknowledgements} The authors are thankful to IMPA and its staff for 
the fine research ambient. A.A and C.M. would like to acknowledge 
Viviane Baladi for the invitation to visit the Institute Henri Poincar\'e in 
May-June 2005, where a large part of the bilinear estimates for the coupling 
terms was done. Also, C.M. is indebted to Terence Tao for some discussions 
about the method of sharp bilinear estimates. A.A. and C.M. were partially 
supported by CNPq-Brazil and A.C. was partially supported by CNPq-Brazil and 
FAPEAL. 

\section{Appendix}\label{s.appendix}  
This appendix collects some well-known results concerning linear and 
multilinear estimates related to the periodic cubic NLS and the periodic KdV, 
and also includes a brief comment about three conserved functionals for the 
NLS-KdV discovered by M. Tsutsumi.  

\subsection{Linear estimates}\label{s.a.linear}We begin with the proof of the 
linear estimates in lemma~\ref{l.linear}. The basic strategy of the argument 
is contained in the work~\cite{CKSTT} of Colliander, Keel, Staffilani, Takaoka 
and Tao. First, we observe that 
$\widehat{\psi U(u_0)}(n,\tau)= \widehat{u_0}(n)\widehat{\psi}(\tau+n^2)$ and 
$\widehat{\psi V(v_0)}(n,\tau)= \widehat{v_0}(n)\widehat{\psi}(\tau-n^3)$. 
Thus, it follows that 
\begin{equation}\label{e.a.linear-0} 
\|\psi(t)U(t)u_0\|_{Z^k}\lesssim\|u_0\|_{H^k} \quad \text{ and } \quad 
\|\psi(t)V(t)v_0\|_{W^s}\lesssim\|v_0\|_{H^s}. 
\end{equation}  
Hence, it remains only to show that 
\begin{equation*} 
\left\|\psi_T(t)\int_0^t U(t-t') F(t') dt'\right\|_{X^k}\lesssim\|F\|_{Z^k} 
\quad \text{ and } \quad \left\|\psi_T(t)\int_0^t V(t-t') 
G(t') dt'\right\|_{Y^s}\lesssim\|G\|_{W^s}. 
\end{equation*}  
Up to a smooth cutoff, we can assume that both $F$ and $G$ are supported on 
$\T\times [-3,3]$. Let $a(t)=\textrm{sgn}(t)\eta(t)$, where $\eta(t)$ is a 
smooth bump function supported on $[-10,10]$ which equals $1$ on $[-5,5]$. The 
identity 
\begin{equation*} 
\chi_{[0,t]}(t') = \frac{1}{2} (a(t')-a(t-t')), 
\end{equation*} 
for $t\in [-2,2]$ and $t'\in [-3,3]$ permits to rewrite 
$\psi_T(t)\int_0^t U(t-t') F(t') dt'$ (resp., 
$\psi_T(t)\int_0^t V(t-t') G(t') dt'$) as a linear combination of 
\begin{equation}\label{e.a.linear-1} 
\psi_T(t)U(t)\int_{\R} a(t')U(-t') F(t') dt' \quad \left(\text{resp., } 
\psi_T(t)V(t)\int_{\R} a(t')V(-t') G(t') dt'\right) 
\end{equation} 
and 
\begin{equation}\label{e.a.linear-2} 
\psi_T(t)\int_{\R} a(t-t')U(t-t') F(t') dt' \quad \left(\text{resp., } 
\psi_T(t)\int_{\R} a(t-t')V(t-t') G(t') dt'\right). 
\end{equation}  

For~(\ref{e.a.linear-1}), we note that by~(\ref{e.a.linear-0}), it suffices to 
prove that 
\begin{equation*} 
\|\int_{\R}a(t')U(-t')F(t') dt'\|_{H^k}\lesssim 
\|F\|_{Z^k} \quad \left(\text{resp., } 
\|\int_{\R}a(t')V(-t')G(t') dt'\|_{H^s}\lesssim \|G\|_{W^s}\right). 
\end{equation*}  
Since the Fourier transform of $\int_{\R}a(t')U(-t')F(t') dt'$ (resp., 
$\int_{\R}a(t')V(-t')G(t') dt'$) at $n$ is 
$\int\widehat{a}(\tau+n^2)\widehat{F}(n,\tau) d\tau$ (resp., 
$\int\widehat{a}(\tau-n^3)\widehat{G}(n,\tau) d\tau$) and 
$|\widehat{a}(\tau)|=O(\la\tau\ra^{-1})$, the desired estimate follows. 
For~(\ref{e.a.linear-2}), we discard the cutoff $\psi_T(t)$ and note that the 
Fourier transform of $\int_{\R} a(t-t')U(t-t') F(t') dt'$ (resp., 
$\int_{\R} a(t-t')V(t-t') G(t') dt'$) evaluated at $(n,\tau)$ is 
$\widehat{a}(\tau+n^2)\widehat{F}(n,\tau)$ (resp., 
$\widehat{a}(\tau-n^3)\widehat{G}(n,\tau)$). Therefore, the decay estimate 
$|\widehat{a}(\tau)| = O(\la\tau\ra^{-1})$ give us the claimed estimate. This 
proves the lemma~\ref{l.linear}.  

\subsection{Trilinear estimates for $(|u|^2 u)$}\label{s.a.trilinear}Next, we 
prove the trilinear estimate in lemma~\ref{l.nls}. The argument is essentially 
contained in the work~\cite{Bourgain} of Bourgain.\footnote{The ``novelty'' 
here is to estimate the contribution of the weighted $L_n^2 L_{\tau}^1$ 
portion of the $Z^k$ norm, although this is not hard, as we are going to see.} 

By definition of $Z^k$, the hypothesis $k\geq 0$ says that it suffices to show 
that 
\begin{equation*} 
\begin{split} 
&\sup_{\|\phi\|_{L_{n,\tau}^2}\leq 1} 
\sum\limits_{n=n_1+n_2-n_3}\int\limits_{\tau=\tau_1+\tau_2+\tau_3} 
\ov{\phi(n,\tau)}\frac{\la n\ra^k}{\la \tau+n^2\ra^{1/2}} 
\widehat{u}(n_1,\tau_1)\widehat{v}(n_2,\tau_2) 
\ov{\widehat{w}(n_3,\tau_3)} \lesssim \\ 
&\|u\|_{X^{k,\frac{3}{8}}} \|v\|_{X^{k,\frac{3}{8}}}\|w\|_{X^{k,\frac{3}{8}}} 
\end{split} 
\end{equation*} 
and 
\begin{equation*} 
\left\|\frac{\la n\ra^k}{\la\tau+n^2\ra} 
\widehat{uv\ov{w}}(n,\tau)\right\|_{L_n^2 L_{\tau}^1}\lesssim 
\|u\|_{X^{k,\frac{3}{8}}}\|v\|_{X^{k,\frac{3}{8}}}\|w\|_{X^{k,\frac{3}{8}}}. 
\end{equation*}  
Observe that $\la n\ra^k\lesssim \max\{\la n_1\ra^k, \la n_2\ra^k, 
\la n_3\ra^k\}$. By symmetry, we can assume that $\la n\ra^k\lesssim 
\la n_1\ra^k$. This reduces matters to show that 
\begin{equation}\label{e.a.trilinear-1} 
\begin{split} 
&\sup_{\|\phi\|_{L_{n,\tau}^2}\leq 1} 
\sum\limits_{n=n_1+n_2-n_3}\int\limits_{\tau=\tau_1+\tau_2+\tau_3} 
\frac{\ov{\phi(n,\tau)}}{\la \tau+n^2\ra^{1/2}} 
\la n_1\ra^k\widehat{u}(n_1,\tau_1)\widehat{v}(n_2,\tau_2) 
\ov{\widehat{w}(n_3,\tau_3)} \lesssim \\ 
&\|u\|_{X^{k,\frac{3}{8}}} \|v\|_{X^{0,\frac{3}{8}}}\|w\|_{X^{0,\frac{3}{8}}} 
\end{split} 
\end{equation} 
and 
\begin{equation}\label{e.a.trilinear-2} 
\begin{split} 
\left\|\sum\limits_{n=n_1+n_2-n_3}\int_{\tau=\tau_1+\tau_2+\tau_3} 
\frac{1}{\la\tau+n^2\ra}\widehat{u}(n_1,\tau_1)\widehat{v}(n_2,\tau_2) 
\ov{\widehat{w}(n_3,\tau_3)}\right\|_{L_n^2 L_{\tau}^1}\lesssim 
\|u\|_{X^{0,\frac{3}{8}}}\|v\|_{X^{0,\frac{3}{8}}}\|w\|_{X^{0,\frac{3}{8}}}. 
\end{split} 
\end{equation}  

First, it is not difficult to see that duality, 
$L_{xt}^4 L_{xt}^4 L_{xt}^4 L_{xt}^4$ H\"older inequality and the 
Bourgain-Strichartz estimate in lemma~\ref{l.Bourgain} (i.e., 
$X^{0,3/8}\subset L^4$) implies~(\ref{e.a.trilinear-1}). Next, consider the 
contribution of~(\ref{e.a.trilinear-2}). By Cauchy-Schwarz in $\tau$, since 
$2(-5/8)<-1$, we need only to prove that 
\begin{equation*} 
\left\|\sum\limits_{n=n_1+n_2-n_3}\int_{\tau=\tau_1+\tau_2+\tau_3} 
\frac{1}{\la\tau+n^2\ra^{\frac{3}{8}}}\widehat{u}(n_1,\tau_1) 
\widehat{v}(n_2,\tau_2) 
\ov{\widehat{w}(n_3,\tau_3)}\right\|_{L_n^2 L_{\tau}^1}\lesssim 
\|u\|_{X^{0,\frac{3}{8}}}\|v\|_{X^{0,\frac{3}{8}}}\|w\|_{X^{0,\frac{3}{8}}}, 
\end{equation*} 
which follows again by duality, $L_{xt}^4 L_{xt}^4 L_{xt}^4 L_{xt}^4$ H\"older 
inequality and the Bourgain-Strichartz estimate. This concludes the proof of 
the lemma~\ref{l.nls}.   

\subsection{Bilinear estimates for $\p_x(v^2)$}\label{s.a.bilinear}Now, we 
present the proof of the bilinear estimate in lemma~\ref{l.kdv}. Since this 
bilinear estimate was used only in the case $s\geq 0$, we will restrict 
ourselves to this specific context (although the proof of the bilinear estimate 
for $-1/2\leq s\leq 0$ is similar). Again, the argument is due to 
Bourgain~\cite{Bourgain} (except for the bound on the weighted 
$L_n^2 L_{\tau}^1$ portion of the $W^s$ norm, which is due to Colliander, 
Keel, Stafillani, Takaoka and Tao~\cite{CKSTT}). By definition of $W^s$, it 
suffices to prove that 
\begin{equation}\label{e.a.bilinear-1} 
\begin{split} 
&\sup_{\|\phi\|_{L_{n,\tau}^2}\leq 1} 
\sum\limits_{n=n_1+n_2}\int_{\tau=\tau_1+\tau_2} 
\frac{|n|\la n\ra^s}{\la\tau-n^3\ra^{1/2}}\widehat{v_1}(n_1,\tau_1) 
\widehat{v_2}(n_2,\tau_2)\ov{\phi(n,\tau)}\lesssim \\ 
&\|v_1\|_{Y^{s,\frac{1}{2}}}\|v_2\|_{Y^{s,\frac{1}{2}-}} + 
\|v_1\|_{Y^{s,\frac{1}{2}-}}\|v_2\|_{Y^{s,\frac{1}{2}}} 
\end{split} 
\end{equation} 
and 
\begin{equation}\label{e.a.bilinear-2} 
\begin{split} 
\left\|\frac{|n|\la n\ra^s}{\la\tau-n^3\ra} 
\widehat{v_1 v_2}(n,\tau)\right\|_{L_n^2 L_{\tau}^1}\lesssim 
\|v_1\|_{Y^{s,\frac{1}{2}}}\|v_2\|_{Y^{s,\frac{1}{2}-}} + 
\|v_1\|_{Y^{s,\frac{1}{2}-}}\|v_2\|_{Y^{s,\frac{1}{2}}}. 
\end{split} 
\end{equation}  

Note that our hypothesis of zero mean implies that $n n_1 n_2\neq 0$. Since 
$$\tau-n^3 = (\tau_1-n_1^3) + (\tau_2-n_2^3) - 3n n_1 n_2,$$ 
we obtain that 
$$\max\{\la\tau-n^3\ra, \la\tau_1-n_1^3\ra, \la\tau_2-n_2^3\ra\}\gtrsim 
|n n_1 n_2|\gtrsim |n|^2.$$  
Also, observe that $s\geq 0$ implies that 
$\la n\ra^s\lesssim \la n_1\ra^s\la n_2\ra^s$.  

First, we deal with~(\ref{e.a.bilinear-1}). To do so, we analyse two cases: 
\begin{itemize} 
\item $\la\tau-n^3\ra=\max\{\la\tau-n^3\ra, \la\tau_1-n_1^3\ra, 
\la\tau_2-n_2^3\ra\}$: in this case, the estimate~(\ref{e.a.bilinear-1}) 
follows from 
\begin{equation*} 
\sup\limits_{\|\phi\|_{L_{n,\tau}^2}\leq 1} 
\sum\limits_{n=n_1+n_2}\int_{\tau=\tau_1+\tau_2} 
\widehat{v_1}\widehat{v_2}\ov{\phi(n,\tau)}\lesssim 
\|v_1\|_{Y^{0,\frac{1}{3}}}\|v_2\|_{Y^{0,\frac{1}{3}}}, 
\end{equation*} 
which is an easy consequence of duality, 
$L_{x,t}^4 L_{x,t}^4 L_{x,t}^2$ H\"older inequality and Bourgain-Strichartz 
estimate in lemma~\ref{l.Bourgain} ($Y^{0,1/3}\subset L^4$). 
\item $\la\tau_j-n_j^3\ra=\max\{\la\tau-n^3\ra, \la\tau_1-n_1^3\ra, 
\la\tau_2-n_2^3\ra\}$ for $j\in\{1,2\}$: in this case, the 
estimate~(\ref{e.a.bilinear-1}) follows from 
\begin{equation*} 
\sup\limits_{\|\phi\|_{L_{n,\tau}^2}\leq 1} 
\sum\limits_{n=n_1+n_2}\int_{\tau=\tau_1+\tau_2} 
\widehat{v_1}(n_1,\tau_1)\widehat{v_2}(n_2,\tau_2) 
\frac{\ov{\phi(n,\tau)}}{\la\tau-n^3\ra^{1/2}} 
\lesssim \|v_1\|_{Y^{0,0}}\|v_2\|_{Y^{0,\frac{1}{2}-}} 
\end{equation*} 
and 
\begin{equation*} 
\sup\limits_{\|\phi\|_{L_{n,\tau}^2}\leq 1} 
\sum\limits_{n=n_1+n_2}\int_{\tau=\tau_1+\tau_2} 
\widehat{v_1}(n_1,\tau_1)\widehat{v_2}(n_2,\tau_2) 
\frac{\ov{\phi(n,\tau)}}{\la\tau-n^3\ra^{1/2}} 
\lesssim \|v_1\|_{Y^{0,\frac{1}{2}-}}\|v_2\|_{Y^{0,0}}, 
\end{equation*} 
which are valid by duality, H\"older and the Bourgain-Strichartz estimate. 
\end{itemize}  

Second, we consider~(\ref{e.a.bilinear-2}). Again, we distinguish two cases: 
\begin{itemize} 
\item $\la\tau_j-n_j^3\ra=\max\{\la\tau-n^3\ra, \la\tau_1-n_1^3\ra, 
\la\tau_2-n_2^3\ra\}$ for $j\in\{1,2\}$: after doing the natural cancelations, 
we see that~(\ref{e.a.bilinear-2}) is a corollary of 
\begin{equation*} 
\begin{split} 
\left\|\la\tau-n^3\ra^{-\frac{2}{3}}\la\tau-n^3\ra^{-\frac{1}{3}} 
\sum\limits_{n=n_1+n_2}\int_{\tau=\tau_1+\tau_2} 
\widehat{v_1}(n_1,\tau_1)\widehat{v_2}(n_2,\tau_2)\right\|_{L_n^2 L_{\tau}^1} 
\lesssim \|v_1\|_{Y^{0,0}}\|v_2\|_{Y^{0,\frac{1}{3}}} 
\end{split} 
\end{equation*} 
and 
\begin{equation*} 
\left\|\la\tau-n^3\ra^{-\frac{2}{3}}\la\tau-n^3\ra^{-\frac{1}{3}} 
\sum\limits_{n=n_1+n_2}\int_{\tau=\tau_1+\tau_2} 
\widehat{v_1}(n_1,\tau_1)\widehat{v_2}(n_2,\tau_2)\right\|_{L_n^2 L_{\tau}^1} 
\lesssim \|v_1\|_{Y^{0,\frac{1}{3}}}\|v_2\|_{Y^{0,0}}. 
\end{equation*}  
Applying Cauchy-Schwarz in $\tau$, since $2(-2/3)<-1$, it suffices to prove 
\begin{equation*} 
\left\|\la\tau-n^3\ra^{-1/3}\sum\limits_{n=n_1+n_2}\int_{\tau=\tau_1+\tau_2} 
\widehat{v_1}(n_1,\tau_1)\widehat{v_2}(n_2,\tau_2)\right\|_{L_{n}^2 L_{\tau}^2} 
\lesssim \|v_1\|_{Y^{0,0}}\|v_2\|_{Y^{0,\frac{1}{3}}} 
\end{equation*} 
and 
\begin{equation*} 
\left\|\la\tau-n^3\ra^{-1/3}\sum\limits_{n=n_1+n_2}\int_{\tau=\tau_1+\tau_2} 
\widehat{v_1}(n_1,\tau_1)\widehat{v_2}(n_2,\tau_2)\right\|_{L_{n}^2 L_{\tau}^2} 
\lesssim \|v_1\|_{Y^{0,\frac{1}{3}}}\|v_2\|_{Y^{0,0}}. 
\end{equation*}  
Rewriting the left-hand sides by duality, using H\"older inequality and 
Bourgain-Strichartz estimate $Y^{0,1/3}\subset L^4$ we finish off this case. 
\item $\la\tau-n^3\ra=\max\{\la\tau-n^3\ra, \la\tau_1-n_1^3\ra, 
\la\tau_2-n_2^3\ra\}$: we subdivide this case into two situations. If 
$\la\tau_j-n_j^3\ra\gtrsim |n n_1 n_2|^{1/100}\gtrsim |n|^{1/50}$ 
for some $j\in\{1,2\}$, we cancel 
$\la\tau_j-n_j^3\ra^{1/6}$ leaving $\la\tau_j-n_j^3\ra^{1/3}$ so that we need 
to show 
\begin{equation*} 
\left\|\la\tau-n^3\ra^{-\frac{1}{2}-} 
\sum\limits_{n=n_1+n_2}\int_{\tau=\tau_1+\tau_2} 
\widehat{v_1}(n_1,\tau_1)\widehat{v_2}(n_2,\tau_2)\right\|_{L_n^2 L_{\tau}^1} 
\lesssim \|v_1\|_{Y^{0,0}}\|v_2\|_{Y^{0,\frac{1}{3}}} 
\end{equation*} 
and 
\begin{equation*} 
\left\|\la\tau-n^3\ra^{-\frac{1}{2}-} 
\sum\limits_{n=n_1+n_2}\int_{\tau=\tau_1+\tau_2} 
\widehat{v_1}(n_1,\tau_1)\widehat{v_2}(n_2,\tau_2)\right\|_{L_n^2 L_{\tau}^1} 
\lesssim \|v_1\|_{Y^{0,\frac{1}{3}}}\|v_2\|_{Y^{0,0}}. 
\end{equation*}  
This is an easy consequence of Cauchy-Schwarz in $\tau$, H\"older inequality 
and Bourgain-Strichartz. If $\la\tau_j-n_j^3\ra\ll |n n_1 n_2|^{1/100}$ for 
$j=1,2$, we observe that 
\begin{equation*} 
\tau-n^3 = -3 n n_1 n_2 + O(\la n n_1 n_2\ra^{1/100}). 
\end{equation*}  
After some cancelations, we need to prove that 
\begin{equation*} 
\left\|\la\tau-n^3\ra^{-1/2}\sum\limits_{n=n_1+n_2}\int_{\tau=\tau_1+\tau_2} 
\widehat{v_1}(n_1,\tau_1)\widehat{v_2}(n_2,\tau_2) \chi_{\Omega(n)}(\tau-n^3) 
\right\|_{L_n^2 L_{\tau}^1} 
\lesssim \|v_1\|_{Y^{0,1/3}}\|v_2\|_{Y^{0,1/3}}, 
\end{equation*} 
where $\Omega(n):=\{\eta\in\R : \eta = -3 n n_1 n_2 + 
O(\la n n_1 n_2\ra^{1/100}) \textrm{ for } n_1,n_2\in\Z \textrm{ with } 
n=n_1+n_2\}$. By Cauchy-Schwarz in $\tau$, we bound the left-hand side by 
\begin{equation*} 
\left\| \left(\int\la\tau-n^3\ra^{-1}\chi_{\Omega(n)}(\tau-n^3) 
d\tau\right)^{1/2} \left\|\sum\limits_{n=n_1+n_2}\int_{\tau=\tau_1+\tau_2} 
\widehat{v_1}(n_1,\tau_1)\widehat{v_2}(n_2,\tau_2) 
\right\|_{L_n^2}\right\|_{L_{\tau}^2}. 
\end{equation*}  
Therefore, it remains only to prove that 
\begin{equation*} 
\left(\int\la\tau-n^3\ra^{-1}\chi_{\Omega(n)}(\tau-n^3) d\tau\right)^{1/2} 
\lesssim 1. 
\end{equation*}  
To estimate the integral on the left-hand side, we need the following lemma 
about the distribution of points in $\Omega(n)$ in a fixed dyadic block: 
\begin{lemma}Fix $n\in\Z-\{0\}$. For $n_1,n_2\in\Z-\{0\}$, we have for all 
dyadic $M\geq 1$ 
\begin{equation*} 
|\{\mu\in\R : |\mu|\sim M, \mu = -3 n n_1 n_2 + O(\la n n_1 n_2\ra^{1/100})\}| 
\lesssim M^{1-\delta}, 
\end{equation*} 
for some $\delta>0$. 
\end{lemma}  
\begin{proof}By symmetry, we may assume $|n_1|\geq |n_2|$. Consider first the 
situation $|n|\geq |n_1|$. Since $\mu = -3 n n_1 n_2 + 
O(\la n n_1 n_2\ra^{1/100})$, we get $|n|\lesssim |\mu|\lesssim |n|^3$ because 
$n_1,n_2\in\Z-\{0\}$ and $|n n_1 n_2|\lesssim |n|^3$. Suppose $\mu\sim M$ and 
$|n|\sim N$. For some $1\leq p\leq 3$, we have $M\sim N^p$. Thus, the 
expression of $\mu$ implies that $|n_1 n_2|\sim M^{1-\frac{1}{p}}$. Observe 
that there are at most $M^{1-\frac{1}{p}}$ multiples of $M^{\frac{1}{p}}$ in 
the dyadic block $\{|\mu|\sim M\}$. Therefore, the set of $\mu$ with the form 
$-3 n n_1 n_2 + O(\la n n_1 n_2\ra^{1/100})$ is the union of 
$M^{1-\frac{1}{p}}$ intervals of size $M^{1/100}$, each of them containing an 
integer multiple of $n$. Then, 
\begin{equation*} 
|\{\mu\in\R: |\mu|\sim M, \mu=-3 n n_1 n_2 + O(\la n n_1 n_2\ra^{1/100})\}| 
\leq M^{1-\frac{1}{p}} M^{1/100}\lesssim M^{3/4}, 
\end{equation*} 
since $1\leq p\leq 3$.  

In the situation $|n|\leq |n_1|$, we have 
$|n_1|\lesssim |\mu|\lesssim |n_1|^3$. So, if $|n_1|\sim N_1$, we obtain 
$M\sim N_1^{p}$ for some $1\leq p\leq 3$.  Thus, we can repeat the previous 
argument. 
\end{proof}  

Using this lemma, it is not hard to prove that 
\begin{equation*} 
\int\la\tau-n^3\ra^{-1}\chi_{\Omega(n)}(\tau-n^3) d\tau\lesssim 1. 
\end{equation*}  
Indeed, we change the variables to rewrite the left-hand side as 
\begin{equation*} 
\int\la\mu\ra^{-1}\chi_{\Omega(n)}(\mu) d\mu. 
\end{equation*}  
Decomposing the domain of integration and using the previous lemma, we have 
\begin{equation*} 
\begin{split} 
&\int\la\mu\ra^{-1}\chi_{\Omega(n)}(\mu) d\mu = \\ 
&= \int\limits_{|\mu|\leq 1}\la\mu\ra^{-1}\chi_{\Omega(n)}(\mu) d\mu + 
\sum\limits_{M\geq 1 \; dyadic}\int\limits_{|\mu|\sim M} 
\la\mu\ra^{-1}\chi_{\Omega(n)}(\mu) d\mu \leq \\ 
&\lesssim 1 + \sum\limits_{M\geq 1 \; dyadic} M^{-1} M^{1-\delta} \lesssim 1. 
\end{split} 
\end{equation*} 
\end{itemize} 
This finishes the proof of the lemma~\ref{l.kdv}.  

\subsection{Three conserved quantities for the NLS-KdV flow}\label{s.a.global}  
In the sequel, we show that the quantities 
\begin{itemize} 
\item $M(t):=\int_{\T} |u(t)|^2 dx$, 
\item $Q(t):=\int_{\T}\left\{\alpha v(t)^2 + 
2\gamma \Im(u(t)\ov{\p_x u(t)}) dx\right\}$ and 
\item $E(t):=\int_{\T}\left\{\alpha\gamma v(t)|u(t)|^2 -\frac{\alpha}{6}v(t)^3 
+ \frac{\beta\gamma}{2}|u(t)|^4 + \frac{\alpha}{2}|\p_x v(t)|^2 + 
\gamma |\p_x u(t)|^2\right\} dx$ 
\end{itemize} 
are conserved by the NLS-KdV flow, as discovered by 
M. Tsutsumi~\cite{MTsutsumi}. By the local well-posedness result in 
theorem~(\ref{t.A}), we may assume that $u$ and $v$ are smooth in both $x$ and 
$t$ variables.  First, we consider $M(t)$. Differentiating with respect to $t$, 
we have 
\begin{equation*} \p_t M(t) = \int_{\T} \p_t u \cdot \ov{u} + \int_{\T} u 
\cdot \ov{\p_t u} 
\end{equation*}  
Since the equation~(\ref{e.nls-kdv}) 
implies 
\begin{equation*} 
\p_t u = i\p_x^2 u -i\alpha uv -i \beta|u|^2 u, 
\end{equation*} 
we see that, by integration by parts, 
\begin{equation*} 
\begin{split} 
\int_{\T}\p_t u \cdot \ov{u} &= i\int\p_x^2 u\cdot\ov{u} 
-i\int\alpha u\ov{u}v -i\int\beta|u|^4 \\ &= -\int u\cdot\ov{i\p_x^2 u} - 
\int u(\ov{-i\alpha uv}) - \int u(\ov{-i\beta|u|^2 u}) \\ 
&= -\int_{\T}u\cdot\ov{\p_t u}. 
\end{split} 
\end{equation*}  
Hence, $\p_t M(t)=0$, i.e., $M(t)$ is a conserved quantity.  Second, we analyse 
$Q(t)$. Differentiating with respect to $t$ and using that $v$ is a real-valued 
function, 
\begin{equation*} 
\p_t Q(t) = 2\alpha\int_{\T} \p_t v\cdot v + 
2\gamma\int_{\T}\Im(\p_t u \ov{\p_x u}) + 
2\gamma\int_{\T}\Im(u\ov{\p_x \p_t u}). 
\end{equation*}  
Applying~(\ref{e.nls-kdv}) and using integration by parts, we obtain 
\begin{equation*} 
2\alpha\int_{\T}\p_t v\cdot v = 
2\alpha\gamma\int_{\T}\{-\p_x^3 v 
-\frac{1}{2}\p_x (v^2)+\gamma\p_x(|u|^2)\}\cdot v = 
2\alpha\gamma\int_{\T}\p_x(|u|^2)\cdot v, 
\end{equation*} 
\begin{equation*} 
\int_{\T} \p_t u\p_x\ov{u} = 
i\int_{\T} \p_x^2u \p_x\ov{u} 
-i\alpha\int_{\T}uv\p_x\ov{u}-i\beta\int_{\T}|u|^2 u\p_x\ov{u} 
\end{equation*} 
and 
\begin{equation*} 
\begin{split} 
\int_{\T}u\p_x \p_t \ov{u} &= \int_{\T} u\cdot\p_x\{-i\p_x^2 u + 
i\alpha \ov{u}v +i\beta|u|^2\ov{u}\}\\ 
&= i\alpha\int_{\T} uv\p_x\ov{u} + i\alpha\int_{\T}|u|^2\p_x v 
+ i\beta\int_{\T}|u|^2 u \p_x\ov{u} +i\beta\int_{\T}|u|^2\p_x(|u|^2) \\ 
&= i\alpha\int_{\T} uv\p_x\ov{u} + i\alpha\int_{\T}\p_x(|u|^2) v 
+ i\beta\int_{\T}|u|^2 u \p_x\ov{u}. 
\end{split} 
\end{equation*}  
In particular, 
\begin{equation*} 
\int_{\T} \p_t u\p_x\ov{u}+\int_{\T}u\p_x \p_t \ov{u} 
= i\int_{\T} \p_x^2u \p_x\ov{u}+i\alpha\int_{\T}\p_x(|u|^2) v. 
\end{equation*}  
Since $i\int_{\T}\p_x^2 u \p_x\ov{u} = \ov{i\int_{\T}\p_x\ov{u}\p_x^2 u}$, we 
get 
\begin{equation*} 
\Im(\int_{\T} \p_t u\p_x\ov{u}+\int_{\T}u\p_x \p_t \ov{u}) = 
\alpha\int_{\T}\p_x(|u|^2) v. 
\end{equation*} 
Hence, putting these informations together, we obtain $\p_t Q(t)=0$.  Third, we 
compute $\p_t E(t)$. Writing $E(t)=I-II+III+IV+V$, where 
$I:=\alpha\gamma\int_{\T}|u|^2 v$, $II=\frac{\alpha}{6}\int_{\T}v^3$, 
$III:= \frac{\beta\gamma}{2}\int_{\T} |u|^4$, $IV:= \frac{\alpha}{2} 
\int_{\T}|\p_x v|^2$ and $V:=\gamma\int_{\T}|\p_x u|^2$.  
Using~(\ref{e.nls-kdv}) and integrating by parts, 
\begin{equation*} 
\p_t I = -\alpha\gamma\int_{\T}|u|^2\{\p_x^3 v+\frac{1}{2}\p_x(v^2)\} 
+v\{i\p_x^2 u\cdot\ov{u} - i\p_x^2\ov{u}\cdot u\}, 
\end{equation*}  
\begin{equation*} 
-\p_t II = -\frac{\alpha}{2}\int_{\T} \left\{-v^2\p_x^3 v + 
\gamma v^2\p_x(|u|^2)\right\} = \frac{\alpha}{2}\int_{\T}v^2\p_x^3 v 
-\alpha\gamma\int_{\T}|u|^2\frac{1}{2}\p_x(v^2), 
\end{equation*}  

\begin{equation*} 
\p_t III = \beta\gamma\int_{\T}\{\p_t u |u|^2\ov{u}+\p_t\ov{u} |u|^2 u\} = 
\beta\gamma\int_{\T}\{i\ov{u}|u|^2\p_x^2 u - iu|u|^2\p_x^2\ov{u} \}, 
\end{equation*}  
\begin{equation*} 
\p_t IV = \alpha\int_{\T}\p_x v \cdot\p_x \p_t v = 
-\frac{\alpha}{2}\int_{\T}v^2\p_x^3 v+ 
\alpha\gamma\int_{\T}|u|^2\p_x^3 v, 
\end{equation*}  
\begin{equation*} 
\begin{split} 
\p_t V &= \gamma\int_{\T}\{\p_x\p_t u\cdot\p_x\ov{u} 
+\p_x u\cdot\p_x\p_t\ov{u}\}  \\ 
&= \gamma\int_{\T}\{i\p_x^2 u 
-i\alpha uv +i\beta |u|^2 u\}\p_x^2\ov{u} 
+ \{-i\p_x^2\ov{u} +i\alpha\ov{u}v - 
i\beta |u|^2\ov{u}\} \p_x^2 u \\ 
&= -\alpha\gamma\int_{\T}v\cdot\{i\p_x^2\ov{u}\cdot u 
-i\p_x^2\ov{u}\cdot u\} -\beta\gamma\int_{\T}\{i\ov{u}|u|^2\p_x^2 u - 
iu|u|^2\p_x^2\ov{u}\}. 
\end{split} 
\end{equation*}  
From these expressions, it is not hard to conclude that $\p_t Q(t)=0$.  


 

\begin{thebibliography}{99} 


\bibitem{Benilov}{E. S. Benilov and S. P. Burtsev,} 
{\it To the integrability of the equations describing the Langmuir-wave-ion- 
wave interaction,} 
{Phys. Let., {\bf 98A} (1983), 256--258.}  

\bibitem{Bourgain}{J. Bourgain,} {\it Fourier transform restriction phenomena 
for certain lattice subsets and applications to nonlinear evolution equations,} 
{Geometric and Functional Anal., {\bf 3} (1993), 107--156, 209--262.}  

\bibitem{CKSTT}{J. Colliander, M. Keel, G. Staffilani, H. Takaoka and T. Tao,} 
{\it Sharp global well-posedness for KdV and modified KdV on $\R$ and $\T$,} 
{J. Amer. Math. Soc., {\bf 16} (2003), 705--749.}  

\bibitem{Corcho}{A. J. Corcho, and F. Linares,} 
{\it Well-posedness for the Schr\"odinger - Kortweg-de Vries system,} 
{Preprint (2005).}  

\bibitem{Funakoshi}{M. Funakoshi and M. Oikawa,} 
{\it The resonant interaction between a long internal gravity wave and a 
surface gravity wave packet,} 
{J. Phys. Soc. Japan, {\bf 52} (1983), 1982--1995.}  

\bibitem{Nishikawa}{H. Hojo, H. Ikezi, K. Mima and K. Nishikawa,} 
{\it Coupled nonlinear electron-plasma and ion-acoustic waves,} 
{Phys. Rev. Lett., {\bf 33} (1974), 148--151.}  

\bibitem{Kawahara}{T. Kakutani, T. Kawahara and N. Sugimoto,} 
{\it Nonlinear interaction between short and long capillary-gravity waves,} 
{J. Phys. Soc. japan, {\bf 39 } (1975), 1379--1386.}  

\bibitem{KPV2}{C. E. Kenig, G. Ponce and L. Vega,} 
{\it A bilinear estimate with applications to the KdV equation,} 
{J. Amer. Math. Soc., {\bf 9} (1996), 573--603.}  

\bibitem{Pecher}{H. Pecher,} 
{\it The Cauchy problem for a Schr\"odinger - Kortweg - de Vries system with 
rough data,} 
{Preprint (2005).}  

\bibitem{Yajima-Satsuma}{J. Satsuma and N. Yajima,} 
{\it Soliton solutions in a diatomic lattice system,} 
{Progr. Theor. Phys., {\bf 62} (1979), 370--378.}  

\bibitem{MTsutsumi}{M. Tsutsumi,} 
{\it Well-posedness of the Cauchy problem for a coupled Schr\"odinger-KdV 
equation,} {Math. Sciences Appl., {\bf 2} (1993), 513--528.}  

\end{thebibliography}
\end{document}